\pdfoutput=1
\documentclass[%
pof,
numerical,
amsmath,
amssymb,
reprint,
]{revtex4-1}
\usepackage{hyperref}
\usepackage{booktabs}
\usepackage{graphicx}
\usepackage{subfigure}
\usepackage{float}
\usepackage[export]{adjustbox}
\usepackage{dcolumn}
\usepackage{bm}
\usepackage{cleveref}
\usepackage[abbreviations]{glossaries-extra}
\usepackage[utf8]{inputenc}
\usepackage[T1]{fontenc}
\usepackage{threeparttable}
\usepackage{mathptmx}
\usepackage{siunitx}
\usepackage{etoolbox}
\usepackage{physics}
\usepackage{color}
\usepackage{CJKutf8}
\usepackage{xcolor}
\usepackage{comment}
\usepackage{multirow}
\usepackage{soul}

\soulregister\cite7
\soulregister\ref7
\soulregister\eqref7
\soulregister\pageref7
\makeatletter
\def\@email#1#2{%
 \endgroup
 \patchcmd{\titleblock@produce}
  {\frontmatter@RRAPformat}
{\frontmatter@RRAPformat{\produce@RRAP{*#1\href{mailto:#2}{#2}}}\frontmatter@RRAPformat}
  {}{}
}%
\makeatother
\begin{document}
\preprint{AIP/123-QED}
\title[]{Vectorial lattice Boltzmann solver for compressible inviscid flows with generic equation of state}
\author{S.A.~Hosseini}
\affiliation{Computational Kinetics Group, Department of Mechanical and Process Engineering, ETH Zürich, 8092 Zurich, Switzerland}%
\author{I.V.~Karlin}
\affiliation{Computational Kinetics Group, Department of Mechanical and Process Engineering, ETH Zürich, 8092 Zurich, Switzerland}%
\email{shosseini@ethz.ch, ikarlin@ethz.ch}
\date{\today}
\begin{abstract}
We develop a vectorial lattice Boltzmann model with a space-time adaptive relaxation coefficient, combined with adaptive time-stepping for compressible Euler dynamics with generic equation of state. The model, due to special form of the relaxation coefficient devised here is shown to converge to the Euler limit with second-order accuracy in the absence of shocks under acoustic scaling. The solver is robust and able to properly capture compressible gas dynamics for both ideal and non-ideal equations of state, including the Bethe--Zel'dovich--Thompson regime. This is verified and demonstrated through a variety of configurations of increasing complexity.
\end{abstract}
\maketitle
\section{Introduction}
\label{sec:Introduction}
Compressible flows at high Mach numbers are central to turbomachinery, hypersonic and re-entry aerodynamics, detonation and blast waves, and astrophysical jets, and their accurate, robust numerical treatment has been a driving problem in computational fluid dynamics for over six decades. The finite-volume, shock-capturing paradigm built on Godunov's scheme~\cite{Godunov1959} and the approximate Riemann solvers of Roe~\cite{Roe1981} and van Leer~\cite{vanLeer1979}, together with the high-resolution total-variation-diminishing (TVD) and essentially-non-oscillatory (ENO) reconstructions that followed, remains the workhorse for the compressible Euler equations~\cite{LeVeque2002,Toro2009,ToroVazquezCendon2012}. These methods span a wide Mach-number range with well-understood accuracy and stability but often come with relatively complex, costly and sensitive numerics, due in part to the presence of non-local and non-linear terms in the balance equations.
The lattice Boltzmann method (LBM) offers a qualitatively different, kinetic route to the same target equations: a small set of discrete-velocity populations relax toward a truncated-Maxwellian equilibrium, and the compressible Navier--Stokes or Euler equations emerge as slow dynamics of the system~\cite{ChenDoolen1998,Succi2001}.
Three decades of development have made LBM a mainstay for low-Mach, incompressible, and turbulent flows, prized for its purely local, parallel update rule; extending it to fully compressible, high-speed flows has been a much more recent undertaking and the subject of a large number of publications in recent years \cite{farag2021unified,frapolli2016lattice,frapolli2015entropic,dorschner2018particles,ji2024eulerian,strassle2026consistent}. While the use of formulations such as the double distribution function approach and the introduction of corrections for leading order errors on standard lattices \cite{prasianakis2008lattice,li2007coupled,saadat2021extended,hosseini2020compressibility} have been successful in allowing for compressible simulations in moderately supersonic flows, classical lattice Boltzmann models still face difficulties for high Mach number flows.
A recently rediscovered alternative has been reported as relatively successful in extending these limits: the \emph{vectorial} lattice Boltzmann method (VLBM), in which each discrete-velocity population is itself vector-valued and carries a share of the \emph{full} conserved
vector $\bm W$ directly, rather than a scalar moment of it. The idea traces to the discrete-velocity kinetic relaxation representations of hyperbolic systems introduced by Jin and Xin~\cite{JinXin1995} and developed systematically as Bhatnagar--Gross--Krook (BGK)-type~\cite{BhatnagarGrossKrook1954} kinetic models by Bouchut~\cite{Bouchut1999,Bouchut2004}; Graille~\cite{Graille2014} cast this as a lattice scheme proper, fixing the vector equilibria by \emph{exact} algebraic consistency with the target flux rather than by moment truncation, and Dubois~\cite{Dubois2014} demonstrated it on strongly nonlinear shallow-water waves. The approach has seen a renewed interest in recent years -- rigorous stability analysis~\cite{GuillonHelieHelluy2024} and positivity- and bound-preserving multidimensional extensions~\cite{WissocqLiuAbgrall2025} chief among them -- that motivates revisiting it as a general-purpose compressible solver.
A Chapman--Enskog expansion of the model shows that it recovers the target Euler system only at leading order; the next-order correction is a genuinely Navier--Stokes-type dissipative term whose strength is set by the relaxation rate $\beta$. Held constant in space and time, and not being subject to any scaling with time-step size, this term is not a vanishing discretization artifact, under acoustic scaling it scales the \emph{same order} as a genuine physical viscous stress. As such it is formally only first-order under acoustic scaling, dissipating both normal and shear signals alike, which makes for a very dissipative solver for the Euler system. Furthermore, extension to a viscous compressible flow solver can prove to be difficult as the leading order error scales similarly to the physical Navier--Stokes level contributions.\\
We propose a practical VLBM solver converging to the compressible Euler balance equation for generic equations of state: This is achieved with a Ducros-gated relaxation sensor scaling with time-step size (\S\ref{sec:num-adaptive}) that sets $\beta(\bm x,t)$ from the local flux divergence, throttled to vanish at the rate required to restore the scheme's underlying accuracy in smooth regions, and gated by a dilatation/vorticity
discriminator~\cite{Ducros1999} so that strongly rotational but smooth structures -- a vortex core, a shear layer -- are not mistaken for shocks. In addition a conservative adaptive time-stepping strategy that maintains propagation on lattice links is proposed in (\S\ref{sec:num-dt-adaptive}). Finally, the proposed solver is agnostic to the equation of state (EOS): because the equilibrium populations are fixed by exact algebraic consistency with the flux $\bm Q(\bm W)$ of the target system rather than by moment truncation, no convexity or specific functional form is required of the EOS, and neither of the first two ingredients need be re-derived when the fluid changes (\S\ref{sec:num-vlbm}). Together, these give a formulation that recovers accurate, low-dissipative Euler-level dynamics for compressible flow with a generic EOS. We demonstrate this on the ideal-gas Euler equations and, as a deliberately severe test of the EOS genericity, a van der Waals gas whose fundamental derivative $\Gamma=1+(\rho/c_s)(\partial c_s/\partial\rho)_s$ changes sign across the Bethe--Zel'dovich--Thompson (BZT) regime~\cite{Thompson1971,Cramer1991,Kluwick2004}, where the classical picture of compressions steepening into shocks and expansions spreading into fans partially inverts. The same sensor and time-stepping strategy, with no problem-specific retuning, are used unmodified across every non-trivial benchmark in \S\ref{sec:validation}, from a one-dimensional shock/entropy-wave interaction to a two-dimensional shock persistently coexisting with a vortex, in both ideal-gas and BZT dense-gas media.\\
The remainder of this paper is organized as follows: We begin in \S\ref{sec:balance} by stating the target balance laws for a generic, single-phase EOS, deferring any specific choice of $P(\rho,e)$ to later sections. \S\ref{sec:numerical} then develops the VLBM scheme itself: \S\ref{sec:num-vlbm} constructs the vectorial equilibrium and the discrete evolution equations; and \S\ref{sec:num-adaptive} introduces the Ducros-gated adaptive relaxation sensor used throughout the remainder of the paper. We proceed in \S\ref{sec:validation} to validate the scheme on a sequence of benchmarks of increasing difficulty: an exact steady vortex, one-dimensional Riemann and composite-wave problems for both ideal and
van der Waals gases, genuinely two-dimensional Riemann problems, a converging BZT compression wave, and, finally, an inviscid shock persistently interacting with a vortex in both an ideal gas and a BZT medium.
\section{Balance equations for a generic compressible fluid}
\label{sec:balance}
We target the \emph{inviscid} compressible flow of a single-component fluid. Viscous and heat diffusion effects are neglected throughout, consistent with the high-speed, weakly dissipative regime of interest. The EOS is left \emph{generic} at this stage: all results below hold for an arbitrary single-phase fluid.
\subsection{Conserved variables}
\label{sec:balance-variables}
The flow is described by the density  $\rho(\bm{x},t)>0$, velocity $\bm{u}(\bm{x},t)\in\mathbb{R}^D$, and specific internal energy $e(\bm{x},t)$, in $D$ spatial dimensions. The specific total energy is
\begin{equation}
  E(\bm{x},t) \;=\; e(\bm{x},t) \;+\; \tfrac12\,|\bm{u}(\bm{x},t)|^2,
  \label{eq:bal-total-energy-def}
\end{equation}
i.e.\ internal plus kinetic energy per unit mass. It is convenient to collect the conserved (extensive, per-unit-volume) variables into the vector
\begin{equation}
  \bm{W} \;=\; \big(\rho,\ \rho\bm{u},\ \rho E\big)^{\dagger} \ \in\ \mathbb{R}^{D+2},
  \label{eq:bal-W}
\end{equation}
namely mass density, momentum density $\rho\bm{u}$, and total energy density $\rho E$. Throughout, $\bm{I}$ denotes the $D\times D$ identity tensor and $\otimes$ the tensor (outer) product, $(\bm{u}\otimes\bm{u})_{ij}=u_iu_j$.
\subsection{Balance equations}
\label{sec:balance-integral}
For an arbitrary fixed control volume $\Omega\subset\mathbb{R}^D$ with boundary $\partial\Omega$ and outward unit normal $\bm{n}$, conservation of mass, momentum, and energy for an inviscid, non-heat-conducting fluid in the absence of body forces reads, in integral form,
\begin{align}
  \dv{t}\int_\Omega \rho \,\dd V
    &= -\oint_{\partial\Omega} \rho\,\bm{u}\cdot\bm{n} \,\dd S,
  \label{eq:bal-int-mass}\\[4pt]
  \dv{t}\int_\Omega \rho\bm{u} \,\dd V
    &= -\oint_{\partial\Omega} \big(\rho\bm{u}\otimes\bm{u} + P\bm{I}\big)\cdot\bm{n} \,\dd S,
  \label{eq:bal-int-momentum}\\[4pt]
  \dv{t}\int_\Omega \rho E \,\dd V
    &= -\oint_{\partial\Omega} \big(\rho E + P\big)\bm{u}\cdot\bm{n} \,\dd S,
  \label{eq:bal-int-energy}
\end{align}
expressing that the rate of change of mass, momentum, and energy in $\Omega$ is due entirely to advective transport and pressure work across $\partial\Omega$ -- no diffusive, viscous, or capillary fluxes are present in the inviscid closure considered here. Applying the divergence theorem and using that $\Omega$ is arbitrary localizes Eqs.~\eqref{eq:bal-int-mass}--\eqref{eq:bal-int-energy}, for
sufficiently smooth solutions, to the differential (strong) form
\begin{align}
  \partial_t \rho + \nabla\cdot(\rho\bm{u}) &= 0,
  \label{eq:bal-mass}\\[4pt]
  \partial_t(\rho\bm{u}) + \nabla\cdot\big(\rho\,\bm{u}\otimes\bm{u} + P\bm{I}\big) &= 0,
  \label{eq:bal-momentum}\\[4pt]
  \partial_t(\rho E) + \nabla\cdot\big[(\rho E + P)\,\bm{u}\big] &= 0,
  \label{eq:bal-energy}
\end{align}
the compressible Euler equations. Equations
\eqref{eq:bal-mass}--\eqref{eq:bal-energy} constitute $D+2$ scalar balance laws for the $D+2$ components of $\bm{W}$, and close as a well-posed system once $P$ is expressed as a function of the thermodynamic state -- the EOS.
It is useful, in particular for the numerical treatment developed later, to write Eqs.~\eqref{eq:bal-mass}--\eqref{eq:bal-energy} compactly as a
single system of conservation laws,
\begin{equation}
  \partial_t \bm{W} \;+\; \nabla\cdot \bm{Q}(\bm{W}) \;=\; \bm{0},
  \label{eq:bal-general-form}
\end{equation}
with conserved vector $\bm{W}$ given by \eqref{eq:bal-W} and flux
\begin{equation}
  \bm{Q}(\bm{W}) \;=\;
  \begin{pmatrix}
    \rho\bm{u} \\[2pt]
    \rho\,\bm{u}\otimes\bm{u} + P\bm{I} \\[2pt]
    (\rho E + P)\,\bm{u}
  \end{pmatrix}.
  \label{eq:bal-flux}
\end{equation}
Written this way, $\bm{Q}$ is an \emph{algebraic} function of $\bm{W}$ alone -- no derivatives of $\bm{W}$ appear on the right-hand side -- so
that \eqref{eq:bal-general-form} is a genuine (first-order, purely hyperbolic) system of conservation laws.
\subsection{Equation of state}
\label{sec:balance-eos}
Equations~\eqref{eq:bal-mass}--\eqref{eq:bal-energy} are not closed: the pressure $P$ must be related to the thermodynamic state through an EOS, which we write generically as
\begin{equation}
  P \;=\; P(\rho,e),
  \label{eq:bal-eos-generic}
\end{equation}
with $\rho$ the density and $e$ the specific internal energy, subject to the usual thermodynamic consistency constraint linking $P(\rho,e)$ to a single-valued entropy $s(\rho,e)$ via the first and second law, $T\,\dd s = \dd e + P\,\dd(1/\rho)$,
\begin{equation}
  \left(\pdv{e}{\rho}\right)_{\!T} \;=\;
  \frac{1}{\rho^2}\left[P - T\left(\pdv{P}{T}\right)_{\!\rho}\right].
  \label{eq:bal-eos-consistency}
\end{equation}
Writing $P_\rho:=(\partial P/\partial\rho)_e$ and
$P_e:=(\partial P/\partial e)_\rho$, the (frozen, isentropic) sound
speed is
\begin{equation}
  c_s^2(\rho,e) \;=\; P_\rho \;+\; \frac{P}{\rho^2}\,P_e,
  \label{eq:bal-soundspeed}
\end{equation}
and the system~\eqref{eq:bal-mass}--\eqref{eq:bal-energy} is hyperbolic provided
\begin{equation}
  c_s^2(\rho,e) \;>\; 0
  \label{eq:bal-hyperbolicity}
\end{equation}
throughout the state space visited by the flow.
We use two closures throughout this work. The \emph{ideal gas},
\begin{equation}
  P \;=\; (\gamma-1)\rho e, \qquad c_s^2 \;=\; \gamma(\gamma-1)e,
  \label{eq:bal-ideal-eos}
\end{equation}
with $\gamma$ the heat-capacity ratio; and the \emph{van der Waals} gas with constant specific heat $c_v$,
\begin{align}
  P(\rho,T) &\;=\; \frac{R\rho T}{1-b\rho} \;-\; a\rho^2,
  \label{eq:bal-vdw-eos}\\
  e(\rho,T) &\;=\; c_v T \;-\; a\rho,
  \label{eq:bal-vdw-caloric}
\end{align}
where $a>0$ and $b>0$ are the attraction and covolume parameters and $R$ the specific gas constant ($\gamma:=1+R/c_v$). Unlike the ideal gas, the
van der Waals fluid's fundamental derivative
$\Gamma=1+(\rho/c_s)(\partial c_s/\partial\rho)_s$ can change sign, giving rise to the non-classical (BZT) wave behavior.
\section{Vectorial lattice Boltzmann model}
\label{sec:numerical}
\subsection{Vectorial lattice Boltzmann formulation}
\label{sec:num-vlbm}
We solve the balance laws of \S\ref{sec:balance} with a VLBM. Unlike standard (scalar-population) lattice Boltzmann schemes, whose populations are scalars recovering the conserved variables only as low-order velocity moments of a truncated Maxwellian expansion, each VLBM population $\bm f_i(\bm x,t)\in\mathbb
R^{D+2}$ is itself vector-valued and carries a share of the \emph{full} conserved vector $\bm W$ of \eqref{eq:bal-W} directly. The populations are transported at a small fixed set of discrete velocities $\bm c_i$, $i=1,\dots,Q$ ($Q$ the number of discrete velocities -- not to be confused with the flux $\bm Q(\bm W)$ of \S\ref{sec:balance}, kept upright/non-bold throughout to distinguish the two), and relax towards an equilibrium $\bm f_i^{\rm eq}(\bm W)$ fixed not by a moment expansion but by \emph{exact} algebraic consistency with the target flux $\bm Q(\bm W)$. This construction is a discrete-velocity relaxation (BGK-type) representation of the hyperbolic system~\eqref{eq:bal-general-form}, in the sense introduced by Jin and Xin~\cite{JinXin1995} for scalar and low-dimensional systems and developed systematically, for general systems of conservation laws and multiple discrete speeds, by Bouchut~\cite{Bouchut1999,Bouchut2004}: the equilibrium is exact (never truncated) for \emph{any} admissible EOS, since it is obtained by solving the moment conditions algebraically rather than by expanding a Maxwellian.
\paragraph{Restriction to the $DdQ2^d$ family.} We restrict throughout to the minimal family denoted $DdQ2^d$ in the standard lattice Boltzmann nomenclature ($d,q$ lower-case placeholders in the family name, instantiated as $D1Q2$ and $D2Q4$ below): $Q=2^D$ discrete velocities for $D\le2$ (the only cases considered in this work), arranged as $D$ antipodal pairs, one pair aligned with each Cartesian direction,
\begin{equation}
  \bm c_{\alpha,-} \;=\; -c\,{\bm e}_\alpha, \qquad
  \bm c_{\alpha,+} \;=\; +c\,{\bm e}_\alpha, \qquad \alpha=1,\dots,D,
  \label{eq:num-velocities}
\end{equation}
with ${\bm e}_\alpha$ the unit vector along direction $\alpha$ and $c$ the (single, isotropic) \emph{link speed} common to every direction. Fixing $c=\delta x/\delta t$ makes the transport sub-step of the scheme an exact grid shift, free of numerical dissipation or dispersion by construction; $\delta x$, $\delta t$ are the (uniform) spatial and temporal discretization steps. The equilibrium consistency conditions are
\begin{equation}
  \sum_{i=1}^{Q} \bm f_i^{\rm eq}(\bm W) \;=\; \bm W,
  \qquad
  \sum_{i=1}^{Q} \bm c_i\,\bm f_i^{\rm eq}(\bm W) \;=\; \bm Q(\bm W),
  \label{eq:num-consistency}
\end{equation}
which lead to the following expression,
\begin{equation}
  \bm f_{i}^{\rm eq}(\bm W) \;=\; \frac{\bm W}{2D} \;+\; \frac{\bm c_i \cdot \bm Q(\bm W)}{2c^2},
  \qquad i=1,\dots,Q,
  \label{eq:num-equilibrium}
\end{equation}
which is the entire EOS content of the scheme: $P$ enters \eqref{eq:num-equilibrium} only through $\bm Q(\bm W)$, so no further modification is required to change the EOS -- ideal gas, van der Waals, or otherwise -- once $\bm Q(\bm W)$ is specified.
\paragraph{Discrete evolution.} One time step splits into an exact transport sub-step,
\begin{equation}
  \bm f_i(\bm x, t) \;=\; \bm f_i(\bm x - \delta t\,\bm c_i,\ t-\delta t),
  \label{eq:num-transport}
\end{equation}
and a relaxation sub-step that over-relaxes each population towards the local equilibrium~\eqref{eq:num-equilibrium},
\begin{equation}
\begin{split}
  \bm f_i(\bm x, t) \;&=\; 2\beta\,\bm f_i^{\rm eq}\big(\bm W(\bm x,t)\big) \\
  &\;+\; (1-2\beta)\,\bm f_i(\bm x,t), \qquad i=1,\dots,Q,
\end{split}
  \label{eq:num-relaxation}
\end{equation}
with $\bm W(\bm x,t)=\sum_{i=1}^Q \bm f_i(\bm x,t)$ and single relaxation parameter $\beta$, common to every direction and population.
Because $\sum_{i=1}^Q \bm f_i^{\rm eq}(\bm W)=\bm W$ holds identically and $\sum_{i=1}^Q\bm f_i(\bm x,t)=\bm W(\bm x,t)$ by definition, the zeroth moment of \eqref{eq:num-relaxation} is invariant under relaxation for \emph{any} $\beta$: mass, momentum, and total energy are conserved exactly at the discrete level, independent of the choice of $\beta$ or of the EOS.
\subsection{Hyperbolicity and linear stability}
\label{sec:num-stability}
\paragraph{Subcharacteristic (hyperbolicity) condition.} For the relaxation system to be a valid, entropy-dissipative approximation of the target system~\eqref{eq:bal-general-form} -- rather than an anti-diffusive, linearly unstable one -- the link speed must dominate every characteristic speed of $\bm Q$ reachable by the flow, in every direction independently~\cite{JinXin1995,Bouchut2004},
\begin{equation}
\begin{split}
  c \;&\ge\; \max_{1\le \alpha\le D}\ \max_{1\le p\le D+2}\ \big|a_p^{(\alpha)}(\bm W)\big|, \\
  a_p^{(\alpha)}(\bm W) \;&=\; \mathrm{eig}_p\big(\bm A_\alpha(\bm W)\big),
\end{split}
  \label{eq:num-subcharacteristic}
\end{equation}
where $\bm A_\alpha(\bm W):=\partial\bm Q_\alpha/\partial\bm W$ is the $\alpha$-th directional flux Jacobian and $a_p^{(\alpha)}(\bm W)$ its $D+2$ eigenvalues (real, by hyperbolicity of the target system, \eqref{eq:bal-hyperbolicity}).
Because a single, isotropic $c$ is shared by every direction in the $DdQ2^d$ family, condition \eqref{eq:num-subcharacteristic} is set by the fastest wave encountered in \emph{any} direction over the whole flow -- in practice $c$ is chosen with a safety margin above $\max(|\bm u\cdot{\bm e}_\alpha|+c_s)$ evaluated over the domain, $c_s$ the sound speed \eqref{eq:bal-soundspeed}.
\paragraph{Linear stability in one dimension.} Linearizing the D1Q2 system about a uniform base state, with $\bm A(\bm W)$ the (locally frozen) flux Jacobian, gives a von-Neumann amplification factor whose physical branch satisfies $|z_{\rm phys}(k)|\le1$ for every wavenumber $k$ if and only if
\begin{equation}
  0 \;<\; \beta \;\le\; 1,
  \label{eq:num-omega-range}
\end{equation}
with $\beta=1$ exactly non-dissipative ($|z_{\rm phys}(k)|=1$ for all $k$) and $\beta=1/2$ recovering the classical Lax--Friedrichs numerical viscosity. We restrict $\beta\in(0,1]$ throughout.
\subsection{Equivalent macroscopic equations}
\label{sec:num-ce}
Because the $D$ directional sub-lattices are mutually independent one-dimensional relaxation systems, each obeys its own Chapman--Enskog expansion exactly as for the single-direction $D1Q2$ system (multiple-scales ansatz $\bm f_i=\bm f_i^{(0)}+\delta t\,\bm
f_i^{(1)}+\delta t^2\bm f_i^{(2)}+O(\delta t^3)$, solvability condition $\sum_{i=1}^Q\bm f_{i}^{(n)}=\bm 0$ for $n\ge1$). Summing the $D+2$ independent results reproduces the target system~\eqref{eq:bal-general-form} exactly at leading order, and gives, at the next order,
\begin{widetext}
\begin{equation}
  \;
  \partial_t \bm W + \nabla\cdot \bm Q(\bm W)
  \;=\; \delta t \sum_{\alpha=1}^{D} \partial_\alpha\!\left\{
  \left(\frac{1}{2\beta}-\frac12\right)
  \Big[c^2\bm I - \bm A_\alpha(\bm W)^2\Big]\,\partial_\alpha\bm W
  \right\} .
  \;
  \label{eq:num-ce-result}
\end{equation}
\end{widetext}
Defining the directional diffusivity matrices
\begin{equation}
\begin{split}
  \bm D_\alpha(\bm W) \;&:=\; \left(\frac{1}{2\beta}-\frac12\right)
  \Big[c^2\bm I - \bm A_\alpha(\bm W)^2\Big], \\
  &\qquad \alpha=1,\dots,D,
\end{split}
  \label{eq:num-diffmatrix}
\end{equation}
Eq.~\eqref{eq:num-ce-result} shows the scheme is consistent with the target Euler system~\eqref{eq:bal-general-form} to $O(\delta t)$, the leading correction being a Fickian-type, \emph{direction-split} (not fully tensorial/rotationally-invariant) dissipative term. Each $\bm D_\alpha(\bm W)$ is positive semi-definite whenever the subcharacteristic condition~\eqref{eq:num-subcharacteristic} and $0<\beta\le1$ hold, since its eigenvalues are
$(c^2-a_p^{(\alpha)2})(1/2\beta-1/2)\ge0$; $\bm D_\alpha\equiv\bm 0$ identically at $\beta=1$.
\subsection{Relaxation coefficient and shock sensor}
\label{sec:num-diffusive-scaling}
\label{sec:num-adaptive}
For $\beta\neq1$, as discussed above and shown in detail in Appendix~\ref{app:multiscale}, the leading order deviation under acoustic scaling from the Euler limit appears at the Navier-Stokes level. While this does not negate convergence to the Euler limit, it means that the scheme is first-order accurate and carries a leading order contribution that operates at the same level as a physical Navier-Stokes dissipation. Possible extensions of the solver to the Navier-Stokes-Fourier limit would then not be a given in that case. In addition, even as a \emph{stabilizing} numerical dissipation mechanism, it is not under our control and is not a shock capturing feature. As a remedy to this issue, as demonstrated in the appendix, $\beta$ must have the following features:
\begin{itemize}
    \item Appendix~\ref{app:multiscale} shows that $\beta$ must follow:
        \begin{equation}
            \frac{1}{2\beta(\bm x,t)} \;=\; \frac12 \;+\; \delta t\,g(\bm x,t),
        \label{eq:num-adaptive-scaling}
        \end{equation}
        where $g\ge0$ is a function to be discussed later. It is easily shown that the above form, in the limit of $\delta t\rightarrow0$ converges to $\beta=1$, and it does so scaling first-order in $\delta t$, pushing the residual dissipative contribution to order three and above.
    \item The function $g$ provides the shock capturing nature of the residual dissipation. We build such a sensor below, gated by one further requirement: a measure built solely from the magnitude of $\nabla\cdot\bm Q(\bm W)$ cannot distinguish a genuine compressive discontinuity from smooth but strongly rotational flow, since both produce large flux gradients, so a dilatation/vorticity discriminator in the spirit of Ducros~\emph{et al.}~\cite{Ducros1999} is used to gate the sensor to compression-dominated regions only.
\end{itemize}
Guided by these two requirements, we build the sensor used throughout this work as follows: For each conserved component $k=1,\dots,D+2$, define the local flux divergence
\begin{equation}
  (\nabla\cdot\bm Q)_k(\bm x) \;:=\; \sum_{\alpha=1}^{D} \partial_\alpha\big(\bm Q_\alpha(\bm W)\big)_k,
  \label{eq:num-sensor-divq}
\end{equation}
evaluated by second-order central differences on the same grid used for transport. We convert it into a local, dimensional estimate of the actual jump in $\bm Q_k$ carried across one grid cell,
\begin{equation}
  j_k(\bm x) \;:=\; \delta x\,(\nabla\cdot\bm Q)_k(\bm x),
  \label{eq:num-sensor-jump}
\end{equation}
and normalize \emph{this} against an absolute flux scale, the
domain-and-direction-wide maximum of $\bm Q$ itself,
\begin{equation}
  \hat s_k(\bm x) \;:=\; \frac{j_k(\bm x)}{\displaystyle\max_{\bm x,\alpha}
  \big|Q_{\alpha,k}(\bm x)\big|},
  \label{eq:num-sensor-normalize}
\end{equation}
placing all $D+2$ fields on a comparable, dimensionless footing
regardless of their physical units. Components are then combined into a single scalar activity measure via an $L^2$ average,
\begin{equation}
  s(\bm x) \;:=\; \left(\frac{1}{D+2}\sum_{k=1}^{D+2} \hat s_k(\bm x)^2\right)^{\!1/2},
  \label{eq:num-sensor-activity}
\end{equation}
and smoothed once, $\bar s:=\mathcal S[s]$, by a separable,
one-cell-radius box filter with weights $(1/4,1/2,1/4)$ applied
successively along each Cartesian direction, to suppress single-cell noise before thresholding.
For the second requirement, the same activity measure $\bar s$ is gated by a dilatation/vorticity discriminator built from the primitive velocity field $\bm u$. In two dimensions, with dilatation and (scalar) vorticity
\begin{equation}
  \mathrm{dil}(\bm x) \;:=\; \partial_x u_x + \partial_y u_y,
  \qquad
  \mathrm{vort}(\bm x) \;:=\; \partial_x u_y - \partial_y u_x,
  \label{eq:num-sensor-dilvort}
\end{equation}
(again by central differences), the discriminator is
\begin{equation}
  \theta(\bm x) \;:=\; \frac{\mathrm{dil}(\bm x)^2}{\mathrm{dil}(\bm x)^2 + \mathrm{vort}(\bm x)^2 + \varepsilon},
  \qquad \bar\theta \;:=\; \mathcal S[\theta],
  \label{eq:num-sensor-theta}
\end{equation}
with $\varepsilon>0$ a small regularizer preventing $0/0$ in regions of exactly quiescent flow; $\theta\to1$ in compression-dominated (shock-like) regions and $\theta\to0$ in rotation/shear-dominated (vortical) regions, since it is built entirely from the ratio of the two and is therefore insensitive to the absolute magnitude of either. In one dimension there is no vorticity to discriminate against and $\bar\theta\equiv1$ identically, so the gate has no effect; we refer to this one-dimensional reduction as Sensor A,
\begin{equation}
\begin{split}
  \chi_A(\bm x) \;&=\; \mathrm{clip}\big(C_{\rm sensor}\,\bar s(\bm x),\ 0,\ 1\big), \\
  \beta_A(\bm x) \;&=\; \beta_{\max} - (\beta_{\max}-\beta_{\min})\,\chi_A(\bm x),
\end{split}
  \label{eq:num-sensor-A}
\end{equation}
with $C_{\rm sensor}>0$ a sensitivity constant and $\mathrm{clip}(\cdot,0,1)$ truncation to $[0,1]$, and to the genuinely two-dimensional construction, retaining the Ducros gate, as Sensor B,
\begin{equation}
\begin{split}
  \chi_B(\bm x) \;&=\; \mathrm{clip}\big(C_{\rm sensor}\,\bar s(\bm x)\,\bar\theta(\bm x),\ 0,\ 1\big), \\
  \beta_B(\bm x) \;&=\; \beta_{\max} - (\beta_{\max}-\beta_{\min})\,\chi_B(\bm x),
\end{split}
  \label{eq:num-sensor-B}
\end{equation}
identical to \eqref{eq:num-sensor-A} except for the extra factor $\bar\theta(\bm x)\in[0,1]$, which suppresses $\chi_B$ (relative to $\chi_A$) wherever the local flux activity is due to smooth rotation rather than compression, while leaving it unchanged at genuine shocks/discontinuities, where $\bar\theta\approx1$. Both variants share the grid-consistent normalization~\eqref{eq:num-sensor-jump}--\eqref{eq:num-sensor-normalize} by construction and are used, without further modification, throughout the remainder of this work. In the results of \S\ref{sec:validation}, $\beta_{\min}$ is set equal to the constant $\beta$ that is otherwise stable for the case at hand (so that neither sensor is ever \emph{more} dissipative, at any point, than the corresponding constant-$\beta$ baseline), and $\beta_{\max}=1$.
\subsection{Adaptive time step}
\label{sec:num-dt-adaptive}
The link speed $c$ has so far been treated as fixed for the duration of a run, chosen once from the initial data via the subcharacteristic condition~\eqref{eq:num-subcharacteristic}, $c\ge\max_{\bm x}\max_\alpha
(|u_\alpha(\bm x,0)|+c_s(\bm x,0))$, with $\delta t=\delta x/c$ fixed correspondingly. This choice is only ever as good as the initial data at predicting the fastest signal speed the flow will reach: if the true maximum signal speed grows substantially during the run (e.g.\ a shock that strengthens, or a flow that accelerates from rest), a $c$ fixed at $t=0$ can under-resolve the subcharacteristic margin later on; if it instead falls (e.g.\ a strong initial transient that subsequently decays or leaves the domain), $c$ remains needlessly large for the remainder of
the run, and with it, via~\eqref{eq:num-diffmatrix}, the numerical dissipation $\bm D_\alpha(\bm W)\propto c^2\bm I-\bm A_\alpha(\bm W)^2$ at every point in the domain -- including points where the flow is otherwise perfectly resolved. We therefore allow $c$, and consequently $\delta t$, to be retargeted during the run rather than fixed at $t=0$.
\paragraph{CFL retargeting.} At (a subset of) time steps, the maximum signal speed currently present in the domain,
\begin{equation}
  s_{\max}(t) \;:=\; \max_{\bm x}\ \max_{1\le\alpha\le D}\
  \big(|u_\alpha(\bm x,t)| + c_s(\bm x,t)\big),
  \label{eq:num-dt-smax}
\end{equation}
is recomputed from the current state, and the link speed is reset to
\begin{equation}
  c(t) \;=\; \frac{s_{\max}(t)}{\mathrm{CFL}}, \qquad
  \delta t(t) \;=\; \frac{\delta x}{c(t)},
  \label{eq:num-dt-retarget}
\end{equation}
with $\mathrm{CFL}\in(0,1)$ a fixed target Courant number setting the safety margin above $s_{\max}(t)$ demanded by
\eqref{eq:num-subcharacteristic}: a smaller $\mathrm{CFL}$ gives a larger margin (larger $c$ relative to $s_{\max}$, more robust to under-resolved transients between retargeting events) at the cost of more numerical dissipation via \eqref{eq:num-diffmatrix}; a larger $\mathrm{CFL}$ the converse. Retargeting every step tracks $s_{\max}(t)$ most tightly; retargeting every few steps amortizes the (small) cost of \eqref{eq:num-dt-smax}--\eqref{eq:num-dt-retarget} at the price of a correspondingly relaxed effective safety margin between updates. The
final step of a run is shortened as needed, $\delta t\to t_{\rm end}-t$, to land exactly on the requested end time.
\paragraph{Reprojection at the new link speed.} Because $c$ enters both the transport sub-step~\eqref{eq:num-transport} and the equilibrium closure~\eqref{eq:num-equilibrium}, changing $c$ mid-simulation requires rebuilding the populations $\bm f_{i}$ consistently with the new value -- they cannot simply be carried over unchanged. A naive rebuild that simply re-evaluates the equilibrium formula at the new speed, $\bm f_{\alpha,i}\leftarrow\bm W/(2D)\mp\bm Q_\alpha(\bm W)/(2c_{\rm
new})$, is \emph{not} correct as the
populations generally carry non-equilibrium content that this substitution silently discards, equivalent to an uncontrolled, spurious full-relaxation step imposed at every retargeting event. The correct reprojection instead holds fixed, across the retargeting event, every \emph{independent moment} the $Q=2^D$ populations actually carry
at that instant -- not just their sum $\bm W$ -- and only changes how those moments are re-expressed as populations at the new link speed.
For the $D2Q4$ lattice used in two dimensions, the four populations $\bm f_{i}$ admit exactly four independent moments: the zeroth moment $\bm W$; the two directional flux-moments
\begin{equation}
  \bm J_\alpha(\bm x,t) \;:=\; \sum_{i=1}^{Q} \bm c_{i\alpha} \bm f_{i}(\bm x,t),
  \qquad \alpha=1,2,
  \label{eq:num-dt-Jalpha}
\end{equation}
which reduce to $\bm J_\alpha=\bm Q_\alpha(\bm W)$ only at equilibrium; and one further, second-order moment,
\begin{equation}
  \bm K(\bm x,t) \;:=\; \sum_{i=1}^{Q} \left( \bm c_{ix}^2-\bm c_{iy}^2\right) \bm f_{i}(\bm x,t),
  \label{eq:num-dt-K}
\end{equation}
which vanishes only at equilibrium ($\bm K^{\rm eq}\equiv\bm0$, since the equilibrium splits $\bm W$ evenly across the two directions by construction). Note $\bm K$ is built
from $c_{ix}^2-c_{iy}^2$ on the four lattice velocities (not $c_{i,x}c_{i,y}$, which vanishes identically on this axis-aligned lattice, every velocity having either $c_{ix}=0$ or $c_{iy}=0$); it is the natural fourth moment completing $\{\bm W,\bm J_1,\bm J_2\}$ into an invertible set for this lattice. Inverting \eqref{eq:num-dt-Jalpha}--\eqref{eq:num-dt-K} together with $\bm W=\sum_{i=1}^{Q} \bm f_{i}$ gives the unique map back to populations at any link speed $c$,
\begin{equation}
 \bm f_i\;=\; \frac{\bm W}{4} + \frac{\bm c_i\cdot\bm J}{2c^2} + \frac{\bm c_{ix}^2-\bm c_{iy}^2}{4c^4}\bm K
  \label{eq:num-dt-invert}
\end{equation}
which reduces exactly to the equilibrium closure~\eqref{eq:num-equilibrium} when $\bm K=\bm0$ and $\bm J_\alpha=\bm Q_\alpha(\bm W)$, and to the single-direction $D1Q2$ formula (no $\bm K$ term, $\alpha=1$ only) when $D=1$.
Retargeting therefore proceeds by evaluating
\eqref{eq:num-dt-Jalpha}--\eqref{eq:num-dt-K} from the \emph{current} populations using the old link speed $c_{\rm old}$, then re-evaluating \eqref{eq:num-dt-invert} at $c_{\rm new}$ with these same values of
$\bm J_1,\bm J_2,\bm K$ held fixed.
\section{Numerical applications}
\label{sec:validation}
We validate the scheme of \S\ref{sec:num-vlbm} on a sequence of benchmarks of increasing difficulty: exact steady solutions, canonical one-dimensional Riemann problems, composite smooth/discontinuous waves, non-ideal (van der Waals, BZT) shock tubes, and genuinely two-dimensional, strongly supersonic flow. We deliberately open with a case containing \emph{no} genuine discontinuity at all, so that the gradient-adaptive relaxation sensor -- introduced specifically to suppress numerical dissipation away from shocks -- can be assessed on its own terms before any shock-capturing demands are placed on it.
\subsection{Static vortex convergence: the Gresho problem}
\label{sec:val-gresho}
\begin{figure}[h]
\centering
\includegraphics[width=0.5\textwidth]{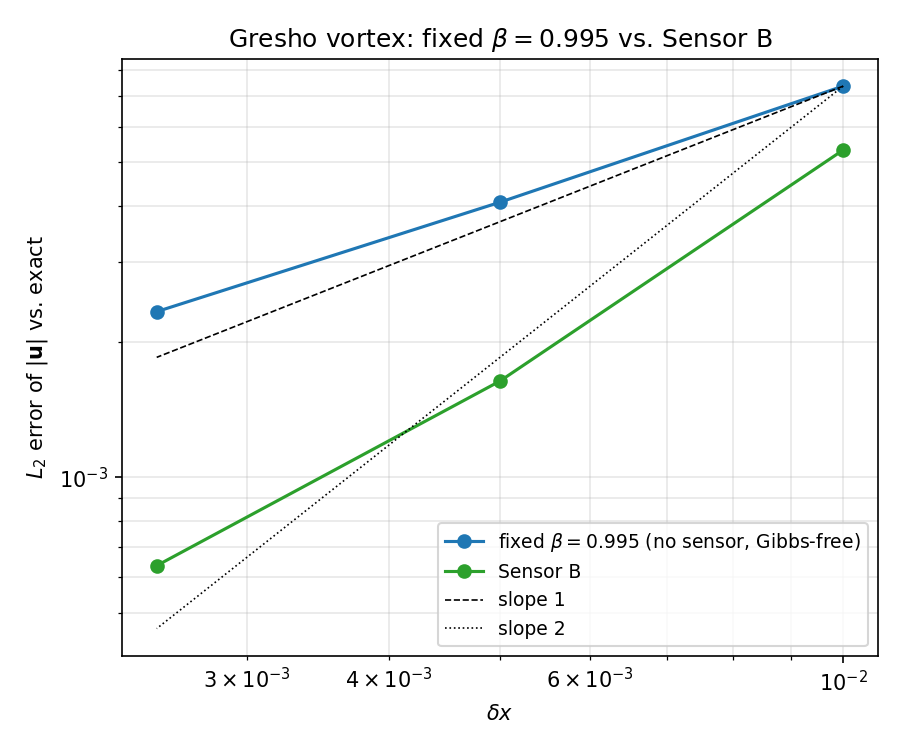}
\caption{Gresho vortex: $L_2$ error of $|\bm u|$ vs.\ $\delta x$ (log--log), Sensor B vs.\ fixed $\beta=0.995$.}
\label{fig:gresho-orders}
\end{figure}

The Gresho vortex~\cite{Gresho1990,LiskaWendroff2003} is the standard test to probe whether the activity of the sensor affects pure shear signals: it is an \emph{exact}, time-independent solution of the compressible Euler equations, so any departure from the initial condition at later time is pure numerical error.\\
On the periodic domain $[0,1]^2$, $\gamma=5/3$, $\rho\equiv1$, the exact steady solution is, in terms of the radius $r=|\bm x-\bm x_0|$ from the domain center $\bm x_0=(0.5,0.5)$,
\begin{widetext}
\begin{equation}
  u_\phi(r) \;=\;
  \begin{cases}
    5r, & r<0.2,\\
    2-5r, & 0.2\le r<0.4,\\
    0, & r\ge0.4,
  \end{cases}
  \qquad
  P(r) \;=\;
  \begin{cases}
    P_0+12.5\,r^2, & r<0.2,\\
    P_0+12.5\,r^2+4\big[1-5r-\ln0.2+\ln r\big], & 0.2\le r<0.4,\\
    P_0-2+4\ln2, & r\ge0.4,
  \end{cases}
  \label{eq:val-gresho-exact}
\end{equation}
\end{widetext}
with $P_0=5$, and $(u_x,u_y)=u_\phi(r)\,(-\sin\phi,\cos\phi)$. The profile is continuous and piecewise-linear in $u_\phi$, with slope discontinuities exactly at $r=0.2$ and $r=0.4$; we run to $t_{\rm end}=0.5$ (several vortex turnover times at the core) on $N\times N$ grids, $N=100,200,400$, and compare the numerical state at $t_{\rm end}$ directly against the (time-independent) exact solution \eqref{eq:val-gresho-exact}.\\
To isolate what the \emph{local} adaptivity of the Sensor buys over simply running less dissipative overall, we compare against the best non-adaptive (spatially and temporally constant $\beta$) scheme achievable for this case: the largest constant $\beta\neq1$ for which the solution remains free of spurious Gibbs-type ringing at the two kinks. Simulations were run using different resolutions. The error scaling with grid size is shown in Fig.~\ref{fig:gresho-orders}.
We can clearly observe that the solver equipped with the sensor exhibits a close-to-second-order scaling of the error with the grid size.
\subsection{Shu--Osher shock/entropy-wave interaction}
\label{sec:val-shuosher}
\begin{figure}[h]
\centering
\includegraphics[width=0.5\textwidth]{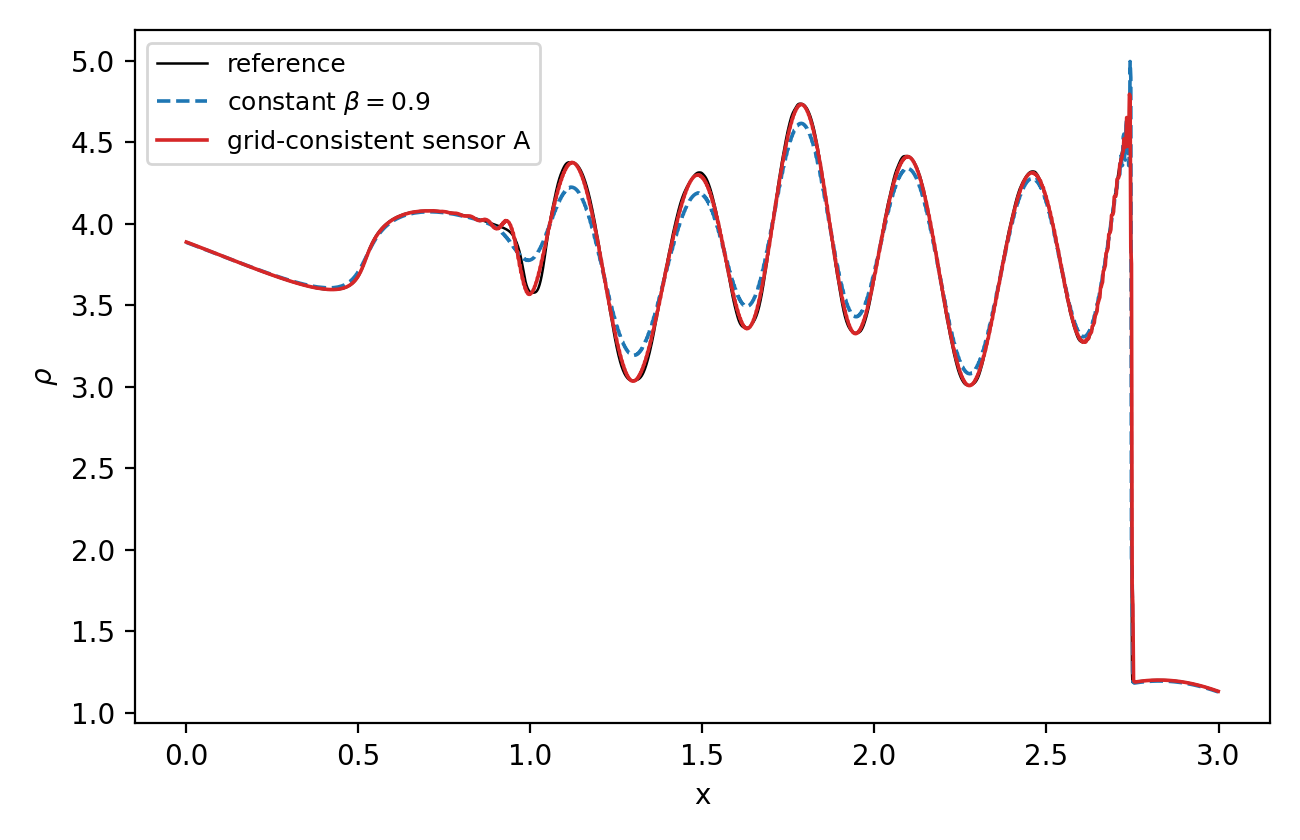}
\caption{Post-shock wave train, zoomed to $x\in[0,3]$, $N=5000$: reference, constant-$\beta$ baseline, and Sensor A.}
\label{fig:shu-osher-zoom}
\end{figure}

\begin{figure}[h]
\centering
\includegraphics[width=0.5\textwidth]{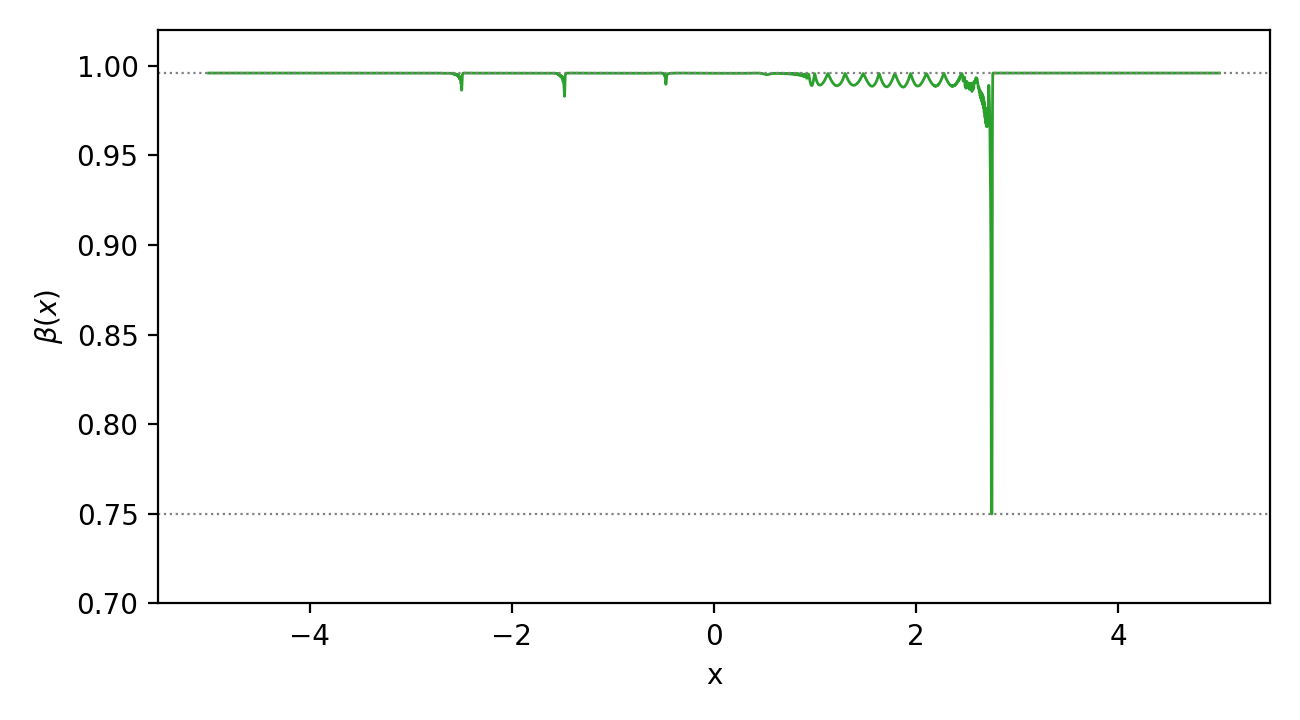}
\caption{Local relaxation parameter $\beta(x)$ for the Shu--Osher problem, $N=5000$.}
\label{fig:shu-osher-beta}
\end{figure}

\begin{figure}[h]
\centering
\includegraphics[width=0.5\textwidth]{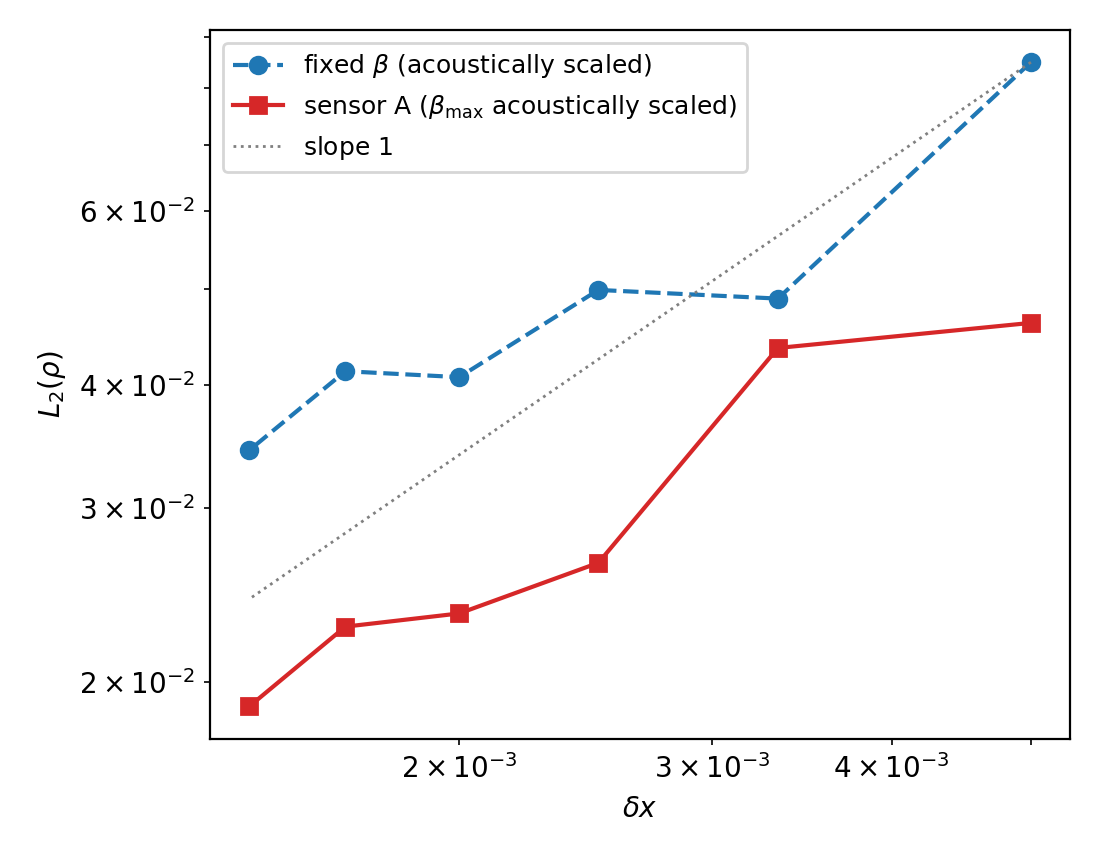}
\caption{Grid convergence, $L_2(\rho)$ error vs.\ $\delta x$, Shu--Osher problem, $N=2000$--$7000$.}
\label{fig:shu-osher-convergence}
\end{figure}

The Shu--Osher problem~\cite{ShuOsher1989} consists of a moving shock driven into a sinusoidally perturbed density field, producing a post-shock region that is genuinely smooth but far from uniform, superposed with a train of compression fronts that steepen as they are advected. A scheme that is too dissipative visibly damps the amplitude of this wave train even though no shock-capturing is required there, while a scheme that is not dissipative enough leaves spurious oscillations at the fronts that do remain.\\
We use the standard Shu--Osher initial condition~\cite{ShuOsher1989} for the ideal-gas Euler equations, $\gamma=1.4$, on $x\in[-5,5]$ with outflow (zero-gradient) boundary conditions at both ends,
\begin{widetext}
\begin{equation}
  (\rho,u,P)(x,0) \;=\;
  \begin{cases}
    (3.857143,\ 2.629369,\ 10.33333), & x<-4,\\[2pt]
    \big(1+0.2\sin(5x),\ 0,\ 1\big), & x\ge-4,
  \end{cases}
  \label{eq:val-shuosher-ic}
\end{equation}
\end{widetext}
run to $t_{\rm end}=1.9$. The left state is the post-shock state of a Mach three shock in the unperturbed gas; as the simulation proceeds this shock sweeps through the sinusoidal density field ahead of it, compressing and advecting it into the wave train.\\
In the absence of a closed-form solution, we validate against an independent high-resolution reference: a second-order Monotonic Upstream-centered Scheme for Conservation Laws (MUSCL) reconstruction (primitive-variable minmod limiter) with a Rusanov numerical flux and explicit second-order Runge--Kutta (RK2) time stepping, run at $N_{\rm ref}=10^4$ grid points, using the same EOS and flux routines as the VLBM solver but an otherwise unrelated numerical method.\\
Two configurations of the D1Q2 VLBM scheme are compared, both using the adaptive time step with target Courant number $\mathrm{CFL}=0.8$, retargeted every step: a constant-relaxation baseline, $\beta\equiv0.9$, and one with our sensor.\\
\begin{figure}[h]
\centering
\includegraphics[width=0.5\textwidth]{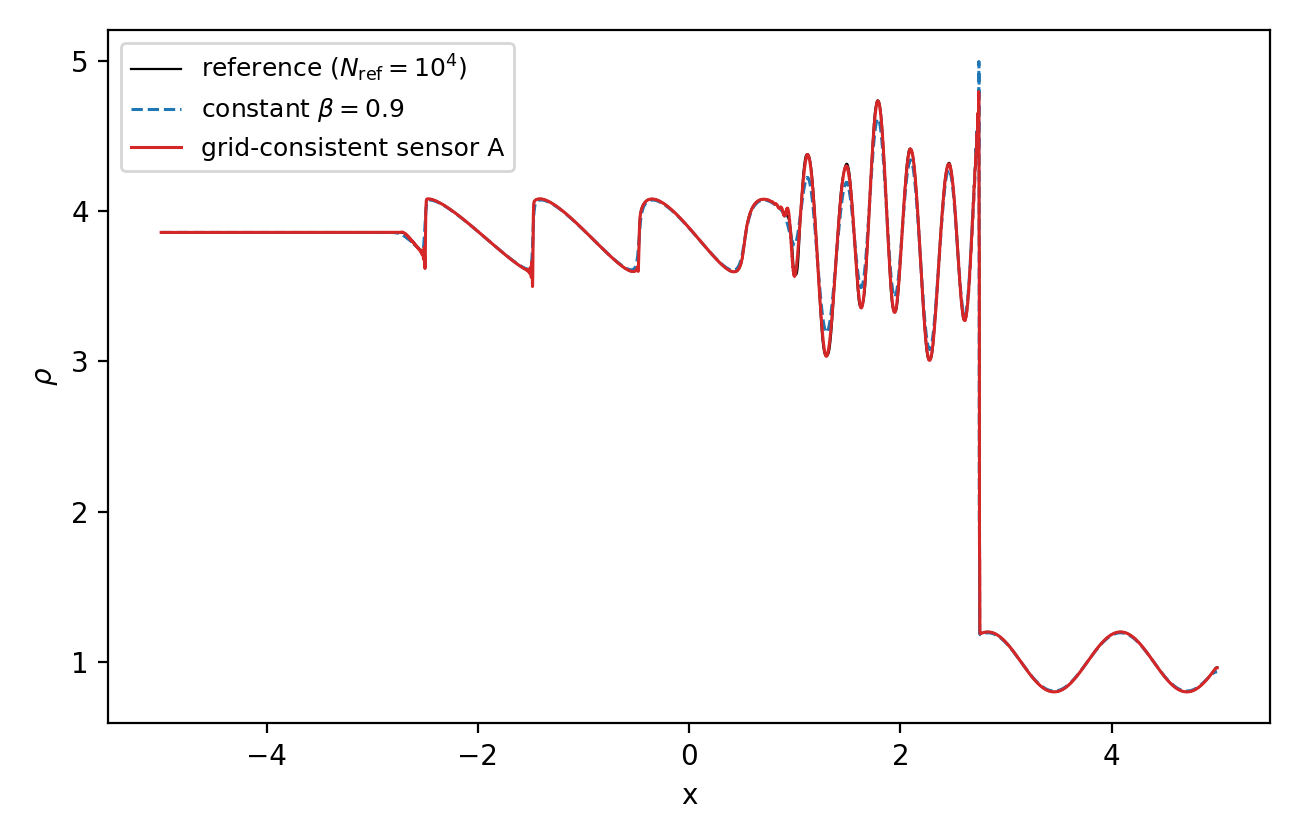}
\caption{Shu--Osher density profile at $t=1.9$, $N=5000$: reference ($N_{\rm ref}=10^4$), constant-$\beta$ baseline, and Sensor A.}
\label{fig:shu-osher-full}
\end{figure}
Figure~\ref{fig:shu-osher-full} shows the full
density profile at $N=5000$ against the reference solution; the shock (at $x\approx2.75$ at $t=1.9$) is captured comparably by both configurations, while the post-shock wave train ($x\gtrsim1$) shows a visible amplitude deficit for the constant-$\beta$ baseline that the sensor-driven run largely removes, most clearly in the zoomed view of Fig.~\ref{fig:shu-osher-zoom}.
Figure~\ref{fig:shu-osher-beta} shows the corresponding local relaxation field $\beta(x)$: it sits at $\beta_{\max}=1$ through essentially the entire smooth wave train, with only shallow, isolated dips of a few grid points' width where the compression fronts are steepest, and drops all the way to the floor $\beta_{\min}=0.75$ only in a narrow region coincident with the main shock. Unlike the previous, relatively normalized sensor, no separate accommodation for the secondary compression fronts immediately behind the shock is needed: the grid-consistent normalization already assigns them an activity level intermediate between the smooth wave train and the fully formed shock.
The convergence plot of Fig.~\ref{fig:shu-osher-convergence} reports $L_2(\rho)$ against the reference at $N=2000,3000,4000,5000,6000,7000$.
While both batches of simulations exhibit a close-to-first-order convergence rate, the solver with the sensor reduces errors by a factor of at least two. The first order convergence is to be expected as at the shock front both solvers become first-order and this error dominates the convergence behavior.
\subsection{Strong shock tube}
\label{sec:val-strong-shock}
\begin{figure*}[t]
\centering
\includegraphics[width=\textwidth]{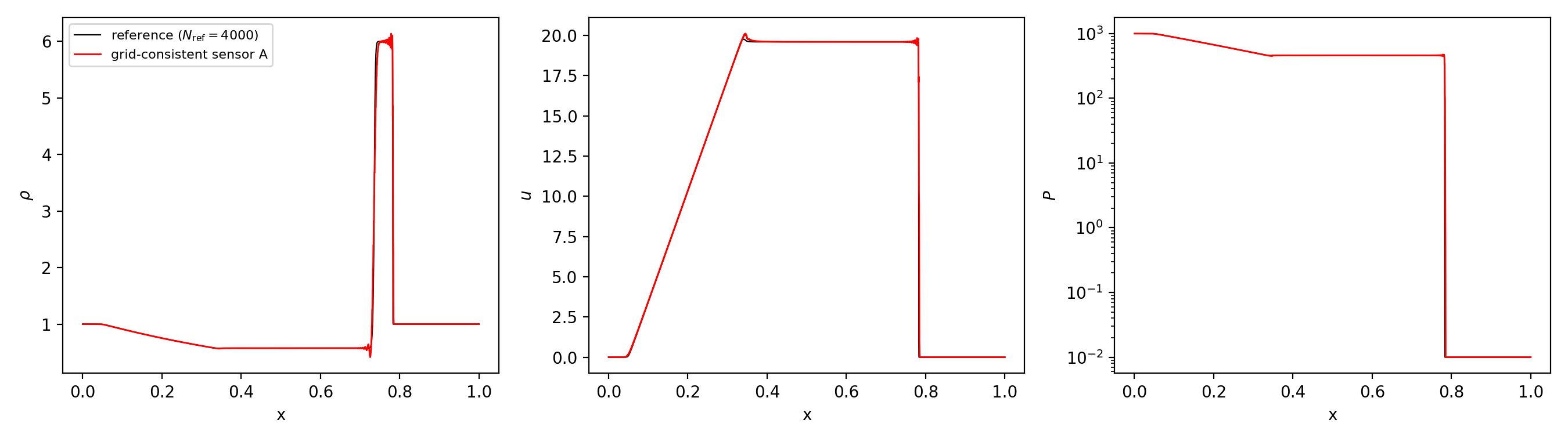}
\caption{Strong shock tube, $N=2000$, $t=0.012$: density, velocity, and pressure against the $N_{\rm ref}=4000$ MUSCL--Rusanov reference.}
\label{fig:strongshock-profiles}
\end{figure*}

\begin{figure}[h]
\centering
\includegraphics[width=0.5\textwidth]{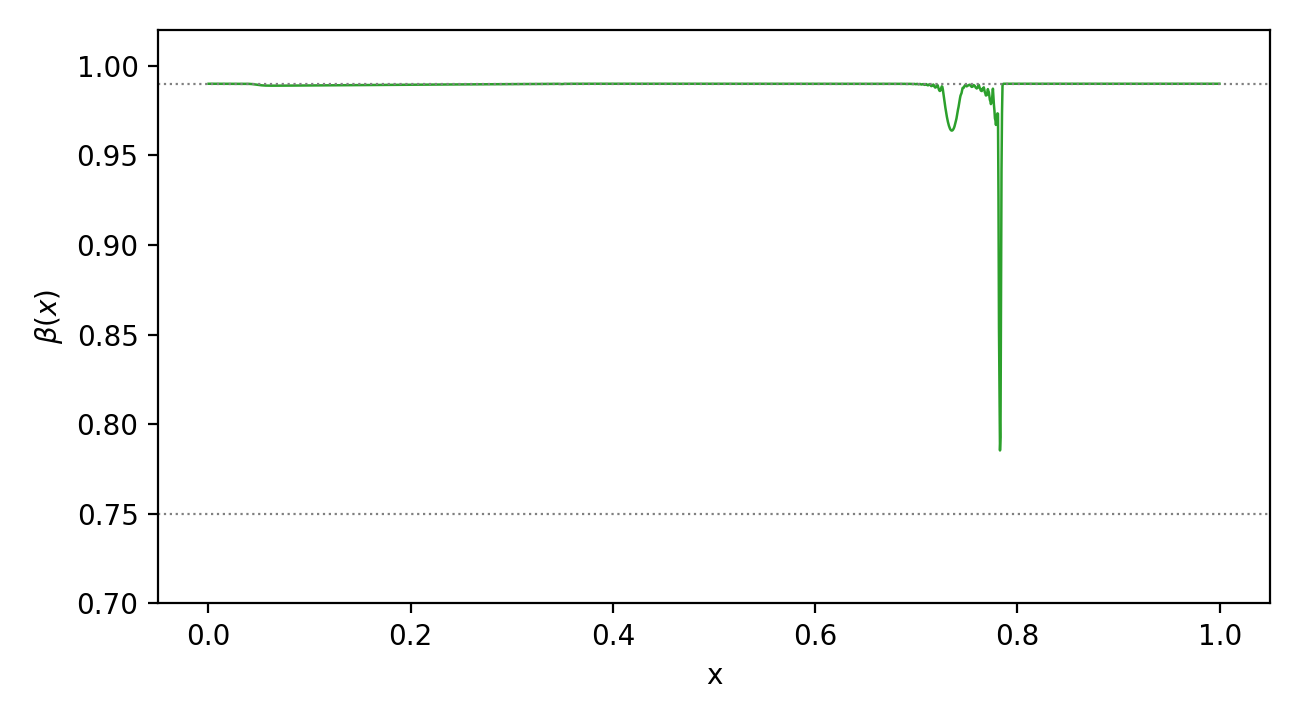}
\caption{Local relaxation parameter $\beta(x)$ for the strong shock tube, $N=2000$.}
\label{fig:strongshock-beta}
\end{figure}

Next, we consider the ``strong shock tube'' problem of Toro and V\'azquez-Cend\'on~\cite{ToroVazquezCendon2012}. This problem is ideal-gas and carries a temperature ratio of $10^5$ between the two initial states and produces a shock of Mach number $\approx 198$ -- a substantially more severe jump than any other benchmark in this paper.\\
The domain is $x\in[0,1]$ with the discontinuity at $x_0=0.5$, ideal gas with $\gamma=1.4$, and initial condition
\begin{equation}
  (\rho,u,P) =
  \begin{cases}
    (1,\,0,\,1000), & x<0.5,\\
    (1,\,0,\,0.01), & x\ge 0.5,
  \end{cases}
  \label{eq:val-strongshock-ic}
\end{equation}
i.e.\ equal densities but a pressure (and hence temperature) ratio of $10^5$. We report the field state at $t=0.012$ with $N=2000$. We validate against an independent MUSCL (primitive-variable minmod reconstruction) plus Rusanov flux, explicit-RK2 reference solver at $N_{\rm ref}=4000$, using the same EOS routines as the VLBM run but an unrelated numerical method.\\
The scheme is run with the adaptive time step ($\mathrm{CFL}=0.8$) and sensor with
\eqref{eq:num-sensor-jump}--\eqref{eq:num-sensor-A}, with
\begin{equation}
  \beta_{\min}=0.75, \qquad \beta_{\max}=1, \qquad C_{\rm sensor}=4.
  \label{eq:val-strongshock-sensor-params}
\end{equation}
Figure~\ref{fig:strongshock-profiles} shows $\rho$, $u$, and $P$ at $t=0.012$ against the reference solution. The expected wave structure is reproduced: a rarefaction fan ($x\lesssim0.35$) connects the initial high-pressure state to an intermediate plateau, a contact discontinuity separates this from a thin, strongly compressed region ($\rho\approx6$, $x\approx0.70$--$0.80$), and a shock of Mach $\approx198$ terminates the disturbance, beyond which the gas remains at its undisturbed initial state. The VLBM profile tracks the independent reference solution closely throughout, with
\begin{equation}
\begin{split}
  L_1(\rho)&=1.39\times10^{-2}, \\
  L_1(u)&=3.33\times10^{-2}, \\
  L_1(P)&=7.18\times10^{-1},
\end{split}
  \label{eq:val-strongshock-l1}
\end{equation}
corresponding to relative errors of $0.2\%$, $0.2\%$, and $0.1\%$ of each field's own maximum value in the reference solution. As with every other benchmark in this work, essentially all of this error is concentrated in the handful of cells spanning the contact and shock themselves (visible in Fig.~\ref{fig:strongshock-profiles} as the rounded rather than perfectly vertical rise into the compressed plateau); a region-by-region breakdown confirms the rarefaction fan and both plateaus contribute negligibly, and the post-shock quiescent region
matches the reference to machine precision.
Figure~\ref{fig:strongshock-beta} shows the corresponding local relaxation field $\beta(x)$: it recovers to near $\beta_{\max}=1$ across the rarefaction fan and both uniform plateaus, and shows two distinct, unfiltered dips in the contact-to-shock region -- a shallow one ($\beta\approx0.96$) at the contact and a much deeper one ($\beta\approx0.79$, close to but not fully at the floor $\beta_{\min}=0.75$) at the shock itself, correctly reflecting that the shock is the stronger of the two features. Both dips are smooth and free of chattering with no post-processing of any kind.
\subsection{Non-ideal (van der Waals) shock tubes: the Argrow benchmarks}
\label{sec:val-argrow-vdw}
\begin{figure}[h]
\centering
\includegraphics[width=0.5\textwidth]{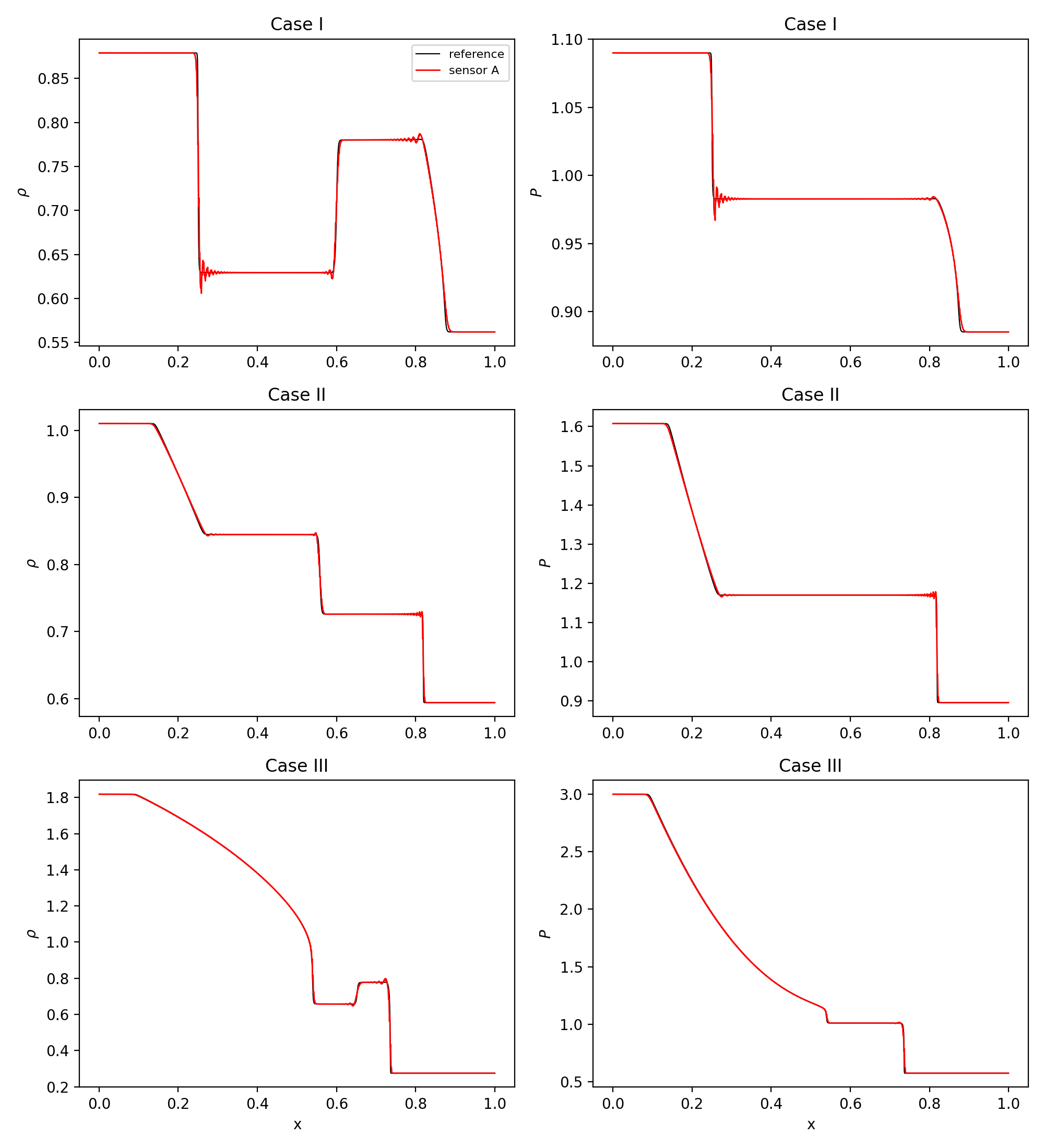}
\caption{$\rho$ (left) and $P$ (right) for the three Argrow van der Waals shock tubes, $N=1000$, against the $N_{\rm ref}=3000$ MUSCL--Rusanov reference.}
\label{fig:argrow-profiles}
\end{figure}

\begin{figure}[h]
\centering
\includegraphics[width=0.5\textwidth]{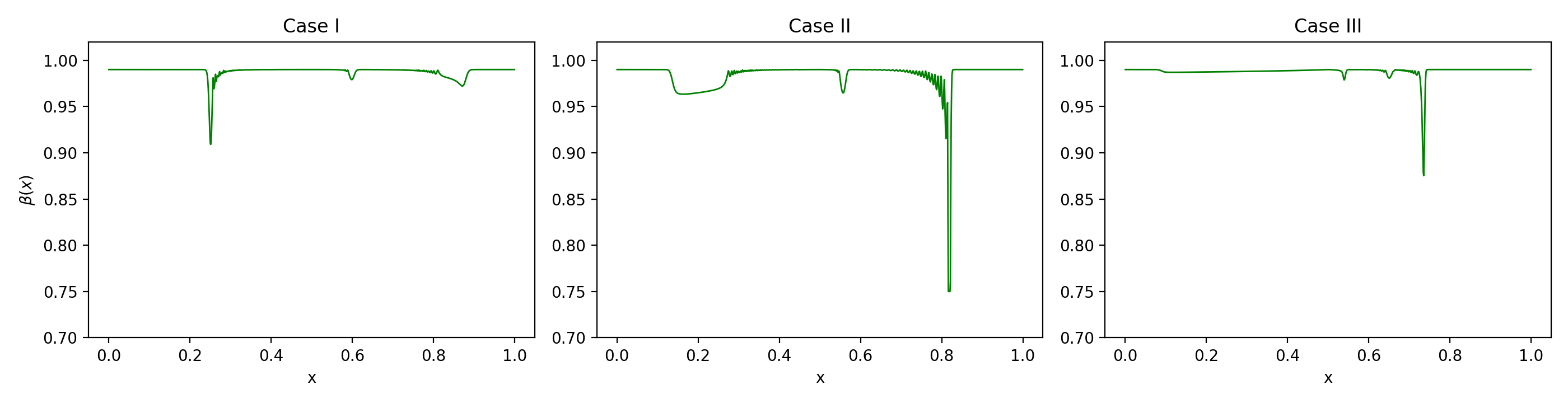}
\caption{Local relaxation parameter $\beta(x)$ for the three Argrow cases, $N=1000$.}
\label{fig:argrow-beta}
\end{figure}

We now leave the ideal-gas closure and validate the scheme on the non-ideal, dense-gas regime that motivates the generic-EOS construction. We use the three van der Waals shock-tube benchmarks of Argrow~\cite{Argrow1996}, subsequently adopted
as a standard validation set for BZT dense-gas solvers~\cite{HosseiniFeinbergKarlin2026}. Unlike the Sod and Shu--Osher problems, these cases are chosen specifically so that the fundamental derivative of gas dynamics,
\begin{equation}
  \Gamma \;=\; 1 + \frac{\rho}{c_s}\frac{\partial c_s}{\partial\rho}\bigg|_s,
  \label{eq:val-argrow-gamma}
\end{equation}
can change sign within the flow: classical gases have $\Gamma>0$ everywhere (compressions steepen into shocks, expansions spread into smooth fans), while a van der Waals gas with sufficiently large specific heat can have $\Gamma<0$ near the critical point, in which case this behavior inverts -- compressions can remain smooth and expansions can steepen into rarefaction shocks. Reproducing this qualitative wave-type inversion, not merely matching a profile, is the substantive test here.\\
All three cases use the van der Waals equation of
state, reduced units ($\rho_c=P_c=1$, i.e.\
$a=3,b=1/3,R=1$), on $x\in[0,1]$ with the discontinuity at $x_0=0.5$ and $u_L=u_R=0$:
\begin{center}
\begin{tabular}{c c c c c c c}
\hline
Case & $R/c_v$ & $P_L$ & $\rho_L$ & $P_R$ & $\rho_R$ & $t_{\rm end}$ \\
\hline
I   & $0.0125$ & $1.09$    & $0.879$ & $0.885$  & $0.562$ & $0.45$ \\
II  & $0.329$  & $1.6077$  & $1.01$  & $0.8957$ & $0.594$ & $0.20$ \\
III & $0.0125$ & $3.00$    & $1.818$ & $0.575$  & $0.275$ & $0.15$ \\
\hline
\end{tabular}
\end{center}
Case II is entirely classical ($\Gamma>0$ throughout); Case I lies close enough to the critical point for the entire sampled state to have $\Gamma<0$; Case III, despite sharing Case I's specific heat, is driven by a much larger pressure jump that carries the flow away from the near-critical $\Gamma<0$ pocket, so it too remains classical, but with a compound-wave structure not seen in Case II. We validate against an independent MUSCL (primitive-variable minmod reconstruction) plus Rusanov flux, explicit-RK2 reference solver at $N_{\rm ref}=3000$, using the same EOS routines as the VLBM run but an unrelated numerical method.\\
All three cases are run with the D1Q2 scheme at $N=1000$, using the adaptive time step ($\mathrm{CFL}=0.8$, retargeted every step).\\
Figure~\ref{fig:argrow-profiles} shows $\rho$ and $P$ for all three cases against the reference solution. The qualitative wave-type inversion is reproduced correctly: Case I shows a smooth
compression fan advancing into the low-pressure state (no jump, purely smooth grid-converged variation across $x\approx0.6$--$0.85$ in $\rho$) together with a sharp rarefaction shock advancing into the high-pressure state ($x\approx0.25$--$0.28$) -- both features direct
consequences of $\Gamma<0$ -- while Case II shows an ordinary shock and rarefaction fan pairing throughout. Case III shows the expected classical, but compound, structure: a leading shock followed by a secondary smooth compression before the trailing contact. In all three cases the VLBM profile with the sensor active tracks the independent reference solution closely, including through the non-classical features
of Case~I, with $L_1(\rho)$ and $L_1(P)$ errors of
$1.31\times10^{-3}$/$4.90\times10^{-4}$ (Case~I),
$6.91\times10^{-4}$/$1.20\times10^{-3}$ (Case~II), and
$1.96\times10^{-3}$/$2.52\times10^{-3}$ (Case~III) at $N=1000$.
Figure~\ref{fig:argrow-beta} shows the corresponding local relaxation field $\beta(x)$ for each case: it dips at every genuine discontinuity -- including the rarefaction shock of Case~I, correctly identified by the flux-divergence activity measure despite it being an expansion wave
rather than a compression -- and recovers to near $\beta_{\max}=1$ across the smooth fans and uniform plateaus in between.
\subsection{Two-dimensional Riemann problems}
\label{sec:val-2d-riemann}
\begin{figure*}[t]
\centering
\includegraphics[width=0.9\textwidth]{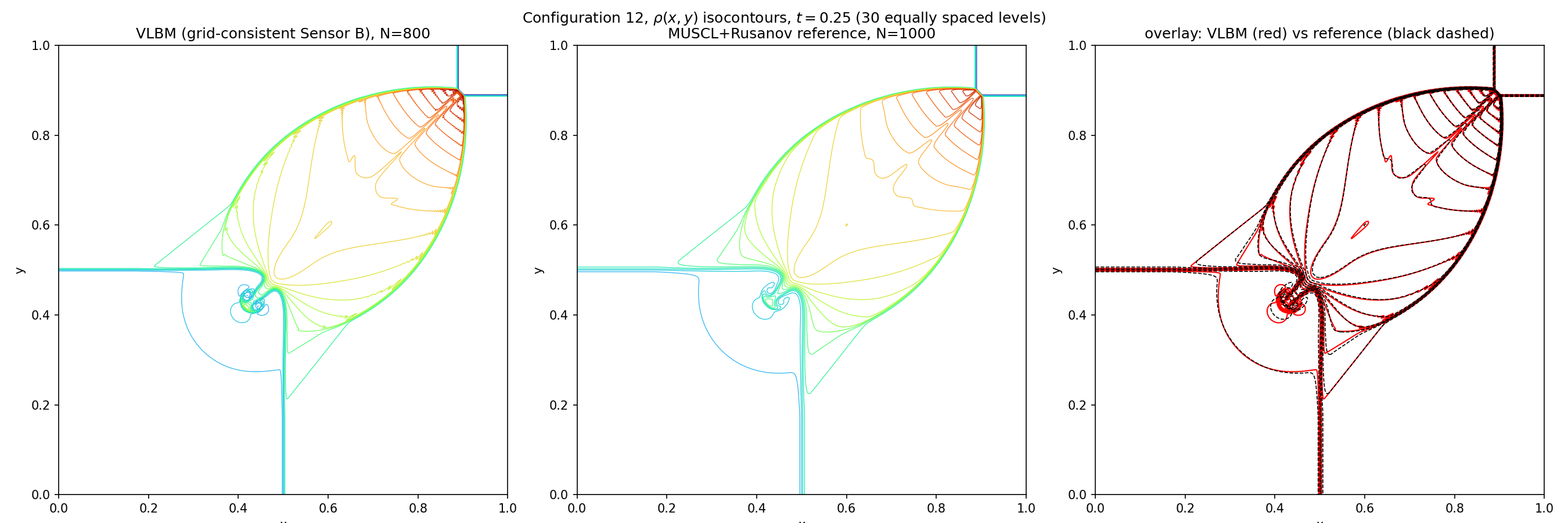}
\caption{Configuration 12, $\rho(x,y)$ at $t=0.25$: 30 equally-spaced density isocontours for the D2Q4 VLBM ($N=800$), the MUSCL+Rusanov reference ($N=1000$), and an overlay.}
\label{fig:riemann12-iso}
\end{figure*}

\begin{figure*}[t]
\centering
\includegraphics[width=0.9\textwidth]{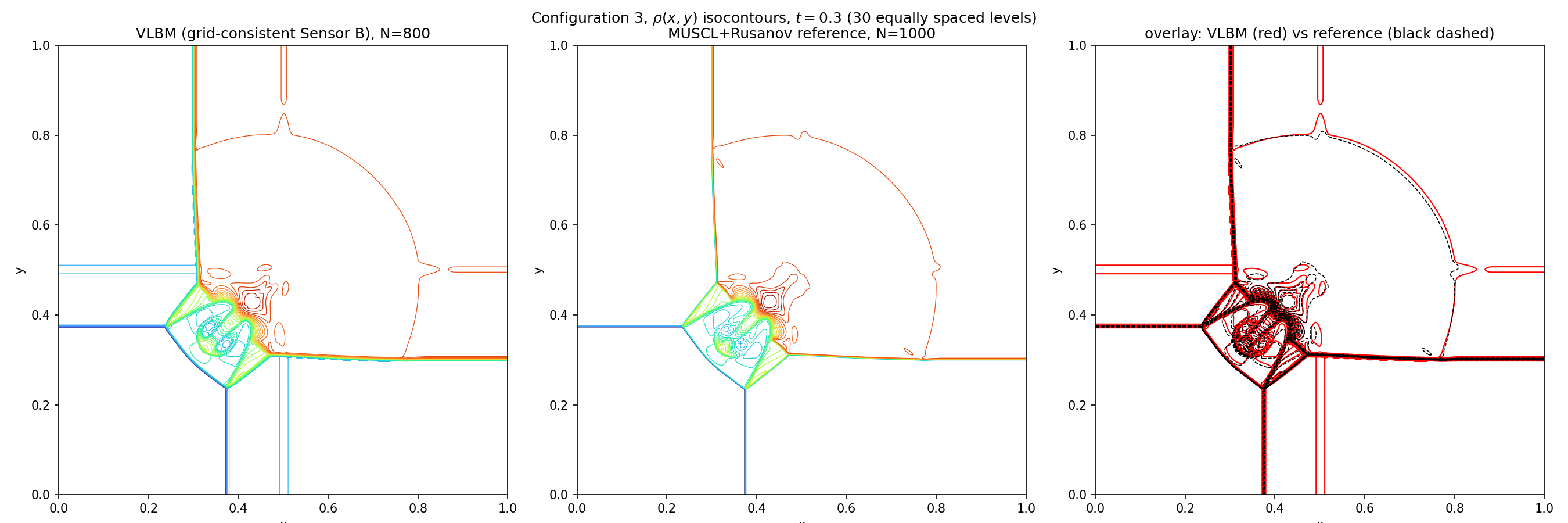}
\caption{Configuration 3, $\rho(x,y)$ at $t=0.3$: 30 equally-spaced density isocontours, same layout as Fig.~\ref{fig:riemann12-iso}.}
\label{fig:riemann3-iso}
\end{figure*}

The benchmarks of \S\S\ref{sec:val-shuosher}--\ref{sec:val-strong-shock} are all one-dimensional; the wave interactions they exercise (shock, contact, rarefaction) are always aligned with the single coordinate axis. To test the D2Q4 scheme on genuinely two-dimensional wave interactions -- oblique shocks, diagonal shock-shock and shock-contact intersections, and the shear-driven roll-up that such intersections can trigger -- we consider two of the classical two-dimensional Riemann problems catalogued by Lax and Liu~\cite{LaxLiu1998} and used extensively since as benchmarks~\cite{KurganovTadmor2002,LiskaWendroff2003}: Configuration 12 and Configuration 3 in the standard numbering.\\
Both problems use the domain $[0,1]\times[0,1]$ with the initial discontinuity along $x=0.5$ and $y=0.5$, ideal gas with $\gamma=1.4$, and constant states in each quadrant
(NE/NW/SW/SE)\@. Configuration 12 has
\begin{equation}
  (\rho,u,v,P) =
  \begin{cases}
    (0.5313,\ 0,\ 0,\ 0.4),      & \text{NE},\\
    (1,\ 0.7276,\ 0,\ 1),        & \text{NW},\\
    (0.8,\ 0,\ 0,\ 1),           & \text{SW},\\
    (1,\ 0,\ 0.7276,\ 1),        & \text{SE},
  \end{cases}
  \label{eq:val-riemann12-ic}
\end{equation}
a mixed shock/contact/contact/rarefaction configuration (edges NE--NW/NE--SE are contacts, NW--SW is a shock, SW--SE is a rarefaction), run to $t_{\rm end}=0.25$. Configuration 3 has
\begin{equation}
  (\rho,u,v,P) =
  \begin{cases}
    (1.5,\ 0,\ 0,\ 1.5),             & \text{NE},\\
    (0.5323,\ 1.206,\ 0,\ 0.3),      & \text{NW},\\
    (0.138,\ 1.206,\ 1.206,\ 0.029), & \text{SW},\\
    (0.5323,\ 0,\ 1.206,\ 0.3),      & \text{SE},
  \end{cases}
  \label{eq:val-riemann3-ic}
\end{equation}
in which \emph{all four} edges are shocks, run to $t_{\rm end}=0.3$. Both stopping times follow the standard convention in the literature~\cite{KurganovTadmor2002} and are short enough that the outermost wave has not yet reached the true domain boundary, so the self-similar quadrant structure near the four corners is unaffected by the boundary treatment; we use outflow (edge-replicated ghost cell) conditions throughout. Both cases are run with the D2Q4 scheme at $N=800\times800$. As in every other benchmark in this work, the VLBM result is cross-validated against an independent solver based on a different numerical mechanism: an unsplit second-order Godunov-type finite-volume scheme (MUSCL reconstruction of the primitive variables with a minmod limiter, Rusanov flux, RK2 time integration), applied directly in two dimensions (both $x$- and $y$-direction fluxes computed from the same reconstructed cell states,
without dimensional splitting), using the same outflow boundary treatment and run to the same $t_{\rm end}$ at $N=1000\times1000$.\\
Figures~\ref{fig:riemann12-iso} and~\ref{fig:riemann3-iso} show 30 equally-spaced density isocontours for the VLBM and reference solutions side by side, together with a direct line-on-line overlay. For Configuration 12, the VLBM reproduces the reference solution's petal-shaped compressed lobe, the curved shock bounding it, and -- most sensitively -- the Kelvin-Helmholtz roll-up spiral at the corner where the two contact discontinuities meet, with the two sets of contours essentially indistinguishable at this scale. Interpolating the VLBM field onto the reference grid gives a relative $L_1$ density error of $0.15\%$ ($L_1=1.52\times10^{-3}$, $L_2=7.16\times10^{-3}$, $L_\infty=2.51\times10^{-1}$, the latter concentrated in the few cells at the spiral core where both solutions have the steepest gradients). For Configuration 3, the four shocks, the diamond-shaped shock-shock interaction region, and the central mushroom-shaped jet produced by the Kelvin-Helmholtz instability at the shock-shock-shock triple point are all reproduced in the correct locations and at the correct scale, with a relative $L_1$ density error of $0.39\%$ ($L_1=3.68\times10^{-3}$, $L_2=1.95\times10^{-2}$,
$L_\infty=4.91\times10^{-1}$)
\subsection{BZT (dense-gas) two-dimensional Riemann problem}
\label{sec:val-bzt-2d}
\begin{figure*}[t]
\centering
\includegraphics[width=\textwidth]{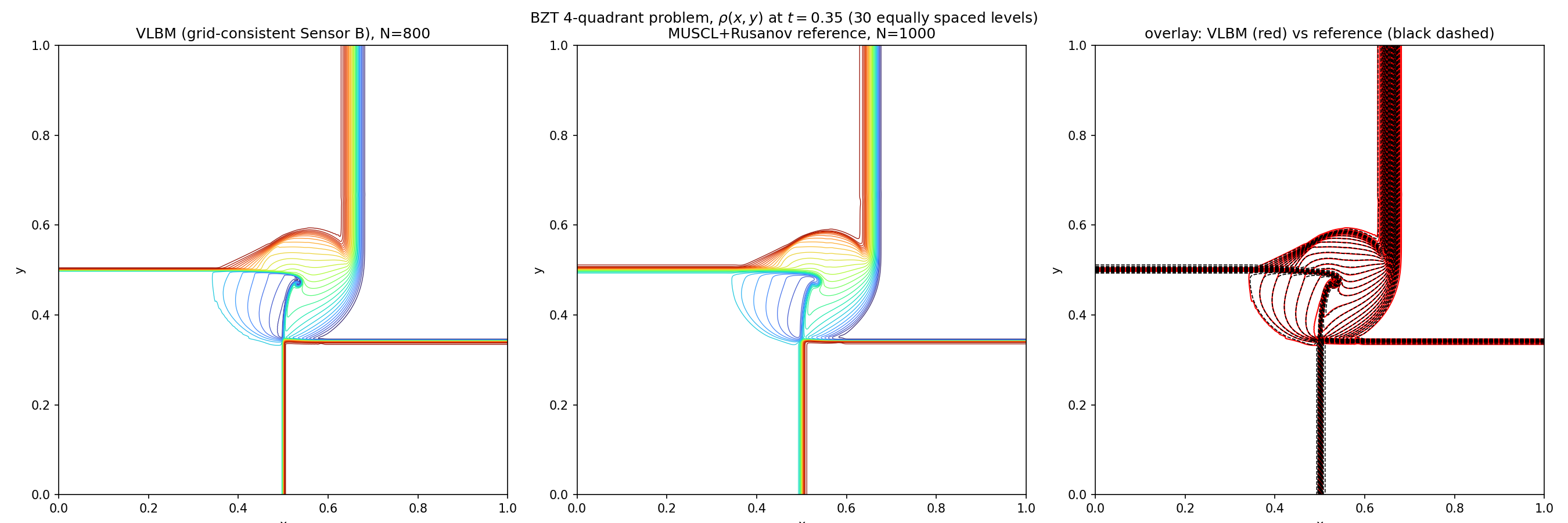}
\caption{BZT 4-quadrant problem, $\rho(x,y)$ at $t=0.35$: 30 equally-spaced density isocontours for the D2Q4 VLBM ($N=800$), the EOS-generic MUSCL+Rusanov reference ($N=1000$), and an overlay.}
\label{fig:bzt-iso}
\end{figure*}

\begin{figure*}[t]
\centering
\includegraphics[width=0.8\textwidth]{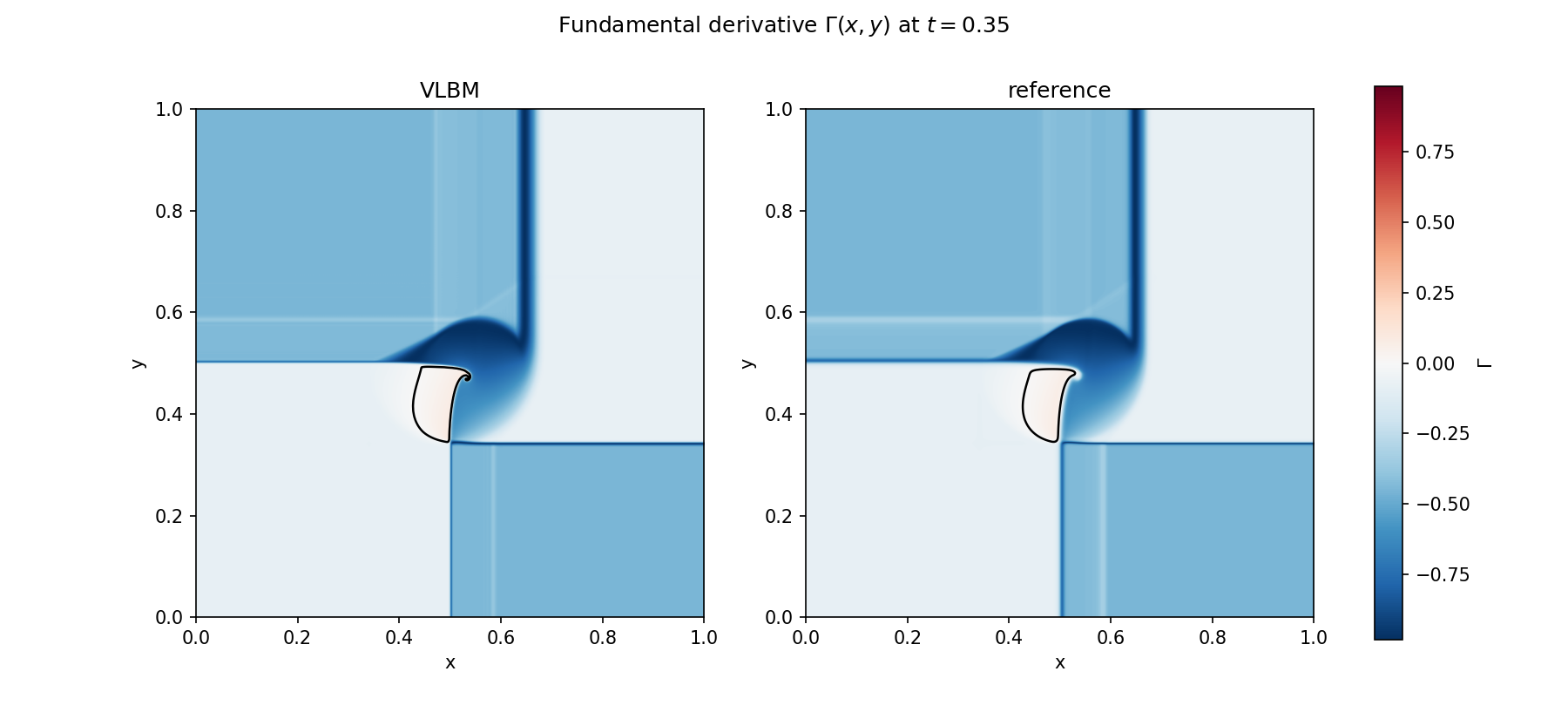}
\caption{Fundamental derivative $\Gamma(x,y)$ at $t=0.35$ for the VLBM (left) and reference (right) solutions; the black contour marks $\Gamma=0$.}
\label{fig:bzt-gamma}
\end{figure*}

The two-dimensional Riemann problems of \S\ref{sec:val-2d-riemann} are both classical ideal-gas configurations, for which the fundamental derivative of gas dynamics
$\Gamma = 1+(\rho/c_s)(\partial c_s/\partial\rho)_s$~\cite{Thompson1971} satisfies $\Gamma=(\gamma+1)/2>0$ everywhere, so compressions always steepen into shocks and expansions always spread into smooth fans. Dense gases with sufficiently large specific heat can instead have $\Gamma<0$ in a neighborhood of the thermodynamic critical point~\cite{Cramer1991,Kluwick2004}. \S\ref{sec:val-argrow-vdw} already validated the scheme's ability to reproduce this inverted wave behaviour in one dimension; here we construct a genuinely two-dimensional BZT Riemann problem with the same 4-quadrant topology as Configuration~12 of \S\ref{sec:val-2d-riemann}, to test the same inversion under oblique, diagonally-interacting waves.\\
We use a van der Waals gas with $a=1.0$, $b=0.3$, $R=1.0$, $\gamma=1+R/c_v=1.02$ (i.e.\ $c_v/R=50$), the same BZT parameter set validated against an independent reference in the
one-dimensional nonclassical shock tube of \S\ref{sec:val-argrow-vdw}; the resulting critical point is $\rho_c=1.111$, $T_c=0.988$, $P_c=0.412$. The four quadrant states,
\begin{widetext}
\begin{equation}
  (\rho,u,v,P) =
  \begin{cases}
    (0.65000,\ 0,\ 0,\ 0.36763),        & \text{NE},\\
    (1.00051,\ 0.15803,\ 0,\ 0.41396),  & \text{NW},\\
    (0.75000,\ 0,\ 0,\ 0.41396),        & \text{SW},\\
    (1.00000,\ 0,\ -0.15747,\ 0.41396), & \text{SE},
  \end{cases}
  \label{eq:val-bzt-ic}
\end{equation}
\end{widetext}
were constructed to reproduce the same wave assignment as Configuration~12 \emph{classically} -- edge NW--SW an exact van der Waals Rankine-Hugoniot shock, edge SE--NE an isentropic simple-wave rarefaction, edges NW--NE and SW--SE contacts (density and tangential velocity free, pressure and normal velocity matched) -- but because $\Gamma<0$ over essentially the whole density range spanned by this problem, the two non-contact edges were independently cross-checked as 1D mini-Riemann problems (VLBM against the MUSCL-Rusanov reference of \S\ref{sec:val-argrow-vdw}) before being assembled into the 2D initial condition, and \emph{both invert}: the classical ``shock'' (NW--SW) actually evolves as a smooth compression fan, and the classical ``rarefaction'' (SE--NE) actually evolves as a rarefaction shock. The domain is $[0,1]\times[0,1]$ with the discontinuity along $x=0.5$, $y=0.5$, run to $t_{\rm end}=0.35$ with outflow boundary conditions (as in \S\ref{sec:val-2d-riemann}, short enough that the outermost wave has not reached the domain boundary), at $N=800\times800$.\\
Figure~\ref{fig:bzt-iso} shows 30 equally-spaced density isocontours for the VLBM and reference solutions. The smooth compression fan (contours spreading apart rather than bunching, upper branch) and the rarefaction shock (contours collapsing onto a single near-vertical front, lower-right branch) are both reproduced, and the overlay shows the two solutions are essentially indistinguishable at this scale, including the small secondary rollup visible where the two inverted waves meet the contact discontinuities. Interpolating the VLBM field onto the reference grid gives a relative $L_1$ density error of $0.12\%$ ($L_1=1.04\times10^{-3}$, $L_2=4.50\times10^{-3}$, $L_\infty=6.80\times10^{-2}$).\\
To confirm directly that this comparison is actually exercising the BZT regime and not merely a van der Waals gas that happens to behave classically, Fig.~\ref{fig:bzt-gamma} shows the fundamental derivative $\Gamma(x,y)$, evaluated pointwise on the final state of both solutions. $\Gamma<0$ over $99.0\%$ of the domain in both the VLBM and the reference solution, confirming that essentially the entire flow field is in the BZT regime. The remaining $\approx1\%$ is a small pocket of classical ($\Gamma>0$) fluid immediately behind the point where the compression fan and the
rarefaction shock intersect, where the local compression pushes the density outside the range originally used to verify $\Gamma<0$ for this EOS; notably, the VLBM and reference solutions agree on this pocket's location, shape, and extent essentially exactly, which is a nontrivial independent cross-check on the sign of $\Gamma$ itself (a third-derivative-level thermodynamic quantity, and correspondingly far more sensitive to numerical error than density) rather than just on the density field.
\subsection{Converging BZT compression wave}
\label{sec:val-bzt-converging}
\begin{figure*}[t]
\centering
\includegraphics[width=0.85\textwidth]{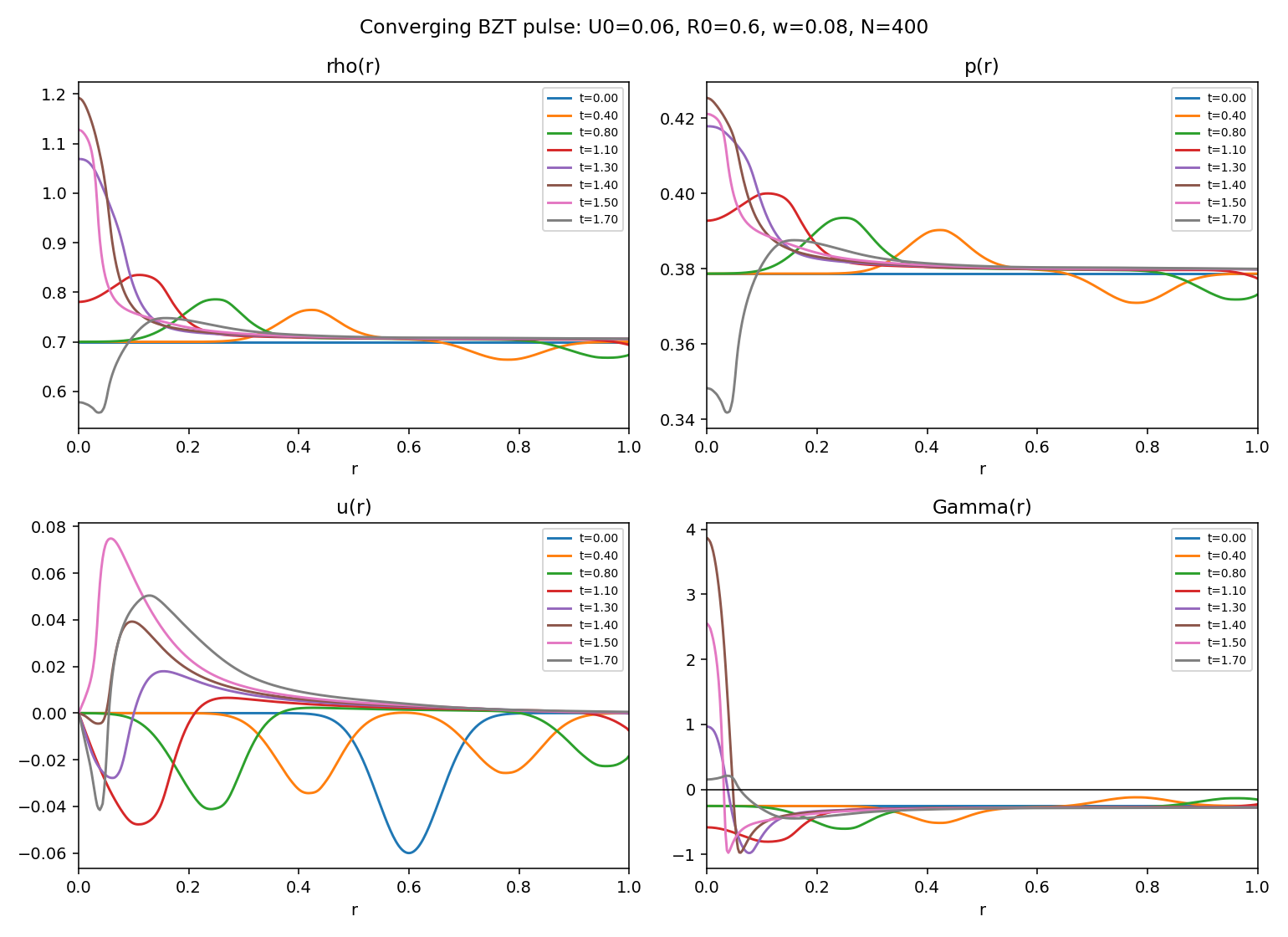}
\caption{1D cylindrical reference: radial profiles of $\rho$, $P$, $u$, $\Gamma$ at several times for the converging BZT pulse ($U_0=0.06$, $R_0=0.6$, $w=0.08$, $N=400$).}
\label{fig:bzt-converging-1d-precheck}
\end{figure*}

\begin{figure*}[t]
\centering
\includegraphics[width=\textwidth]{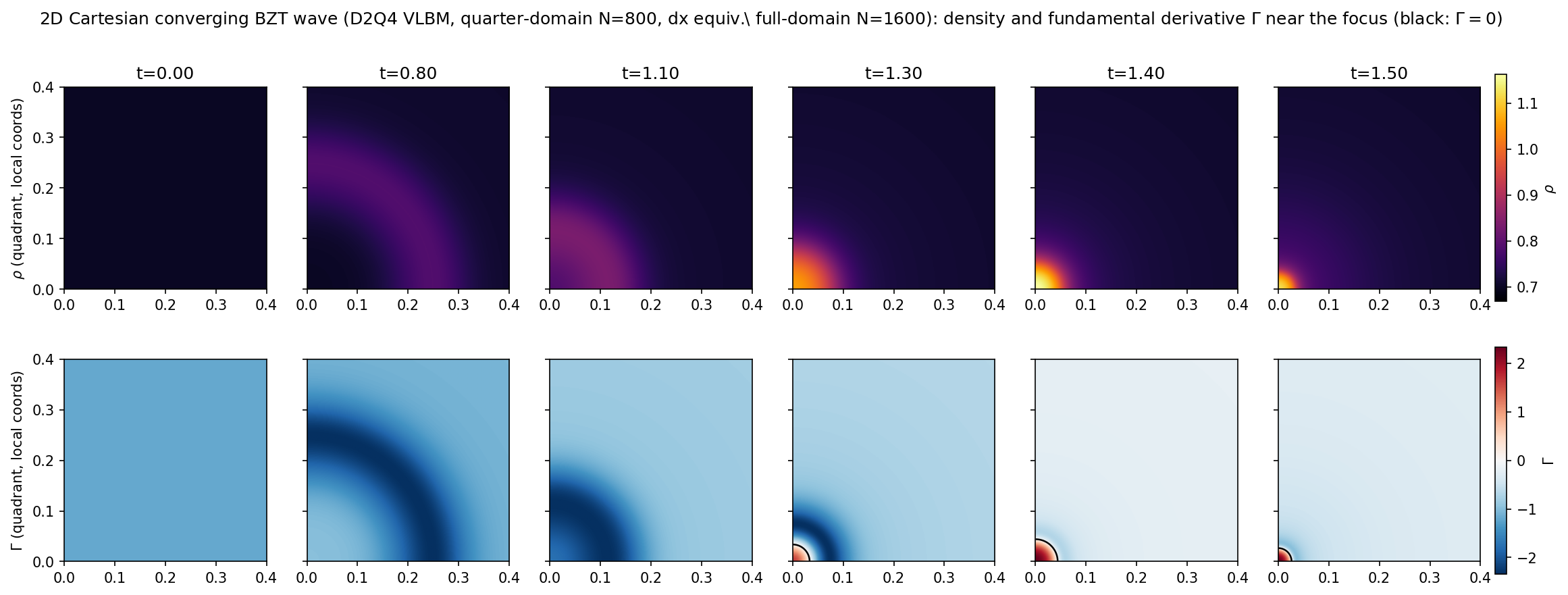}
\caption{2D Cartesian D2Q4 VLBM ($N=400$), converging BZT pulse: density (top) and fundamental derivative $\Gamma$ (bottom) near the focus at several times.}
\label{fig:bzt-converging-2d-captures}
\end{figure*}

\begin{figure*}[t]
\centering
\includegraphics[width=\textwidth]{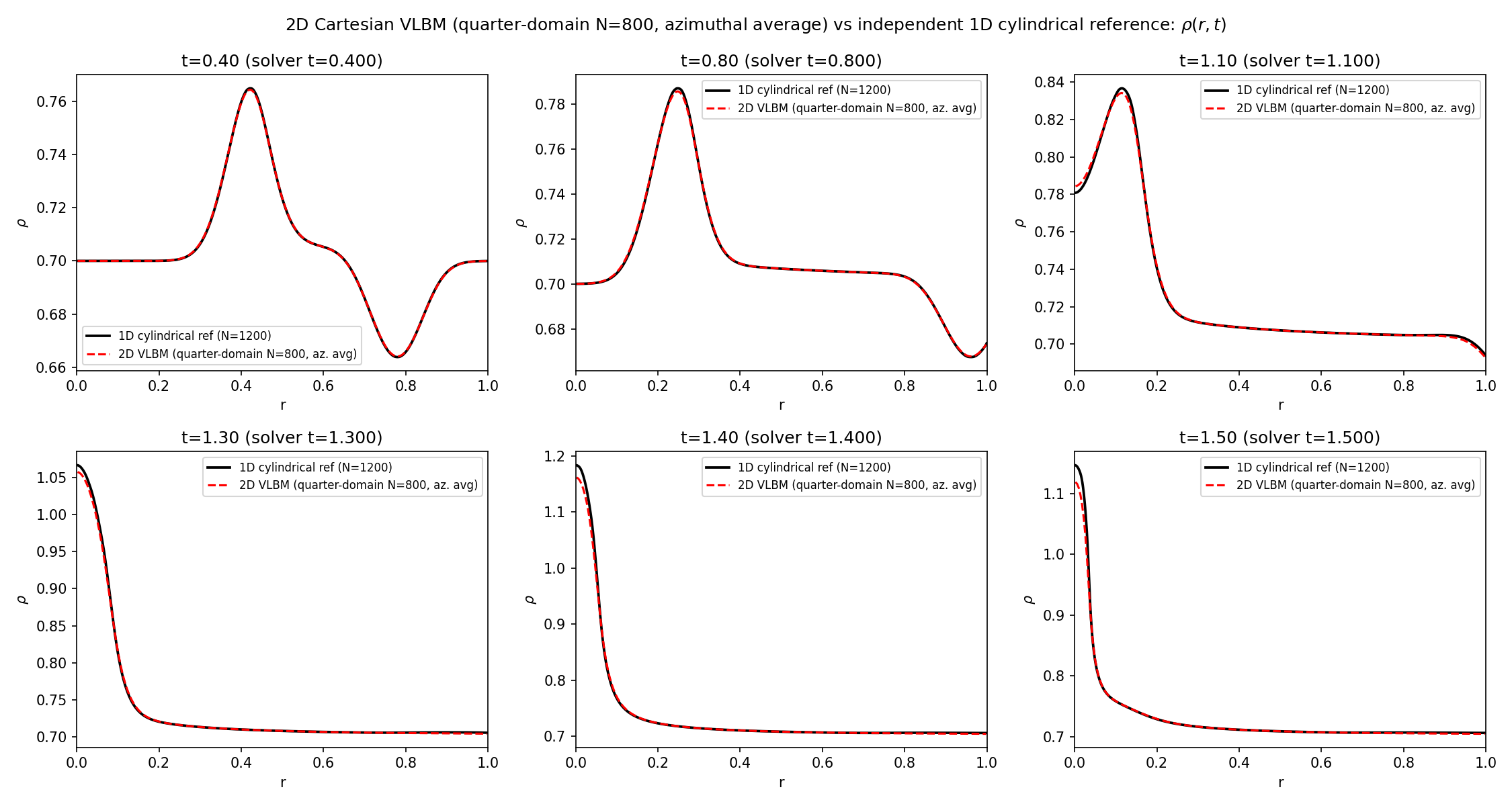}
\caption{Azimuthal average of the 2D Cartesian VLBM density field vs.\ the 1D cylindrical reference ($N=1200$), at matching times.}
\label{fig:bzt-converging-2d-vs-1d}
\end{figure*}

\begin{figure*}[t]
\centering
\includegraphics[width=\textwidth]{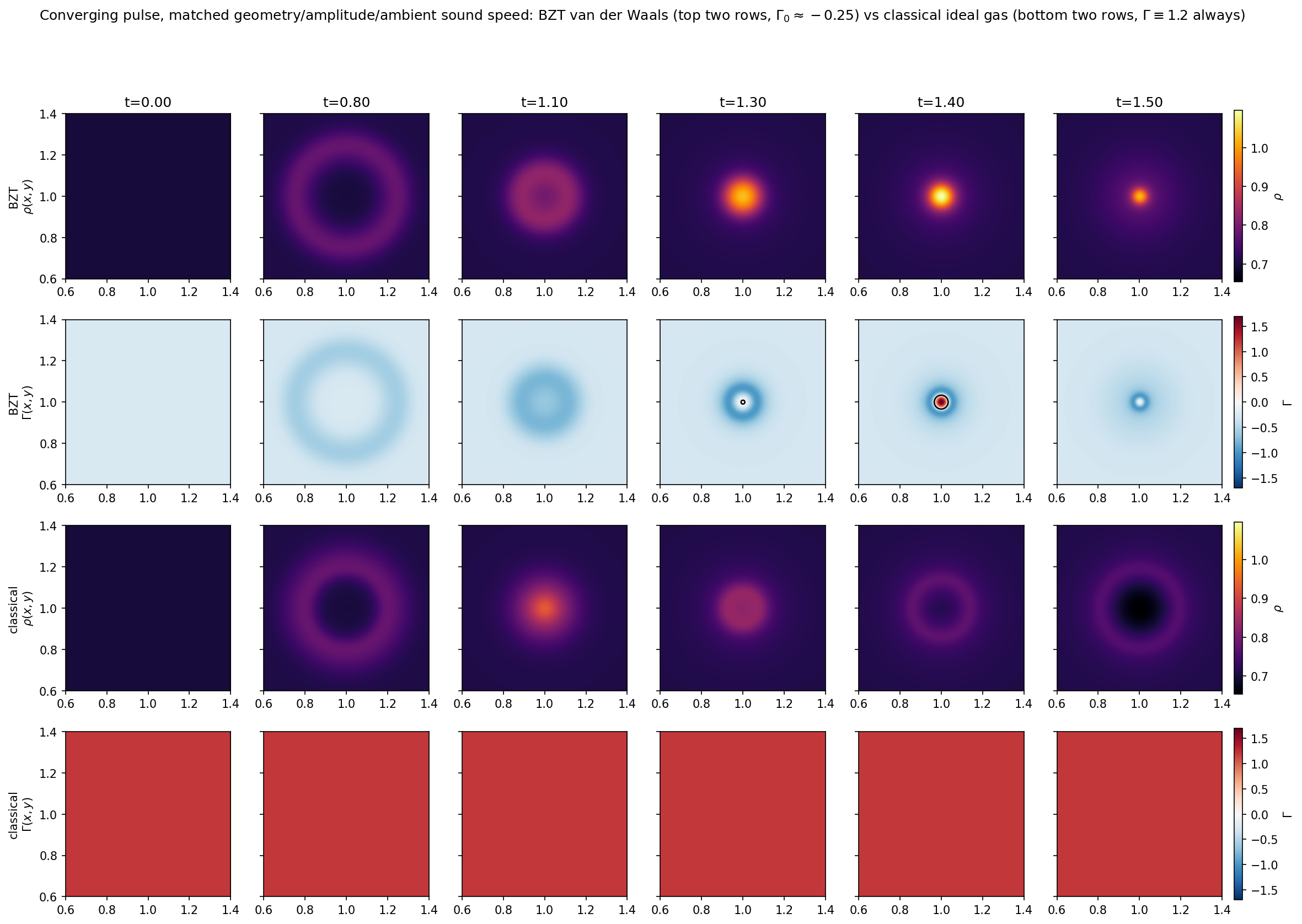}
\caption{Converging pulse with matched geometry, amplitude, and ambient sound speed: BZT van der Waals (top two rows) vs.\ classical ideal gas (bottom two rows).}
\label{fig:bzt-vs-classical-2d}
\end{figure*}

\begin{figure*}[t]
\centering
\includegraphics[width=0.5\textwidth]{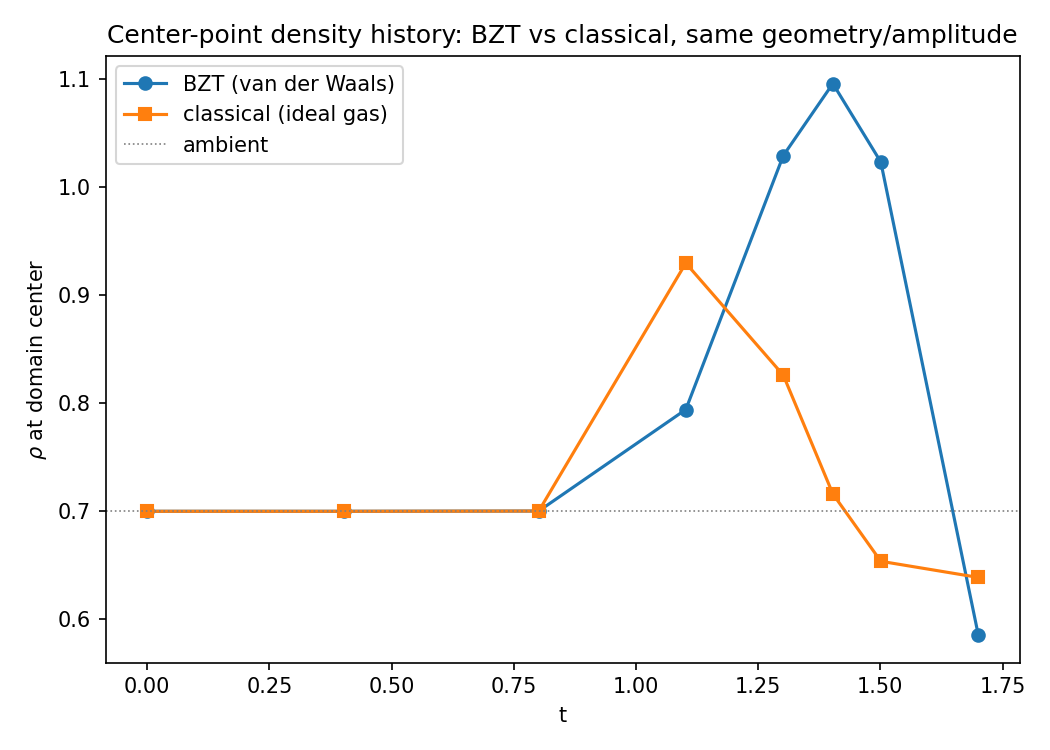}
\caption{Domain-center density history, BZT vs.\ classical, identical geometry/amplitude/ambient sound speed.}
\label{fig:bzt-vs-classical-history}
\end{figure*}

\begin{figure*}[t]
\centering
\includegraphics[width=\textwidth]{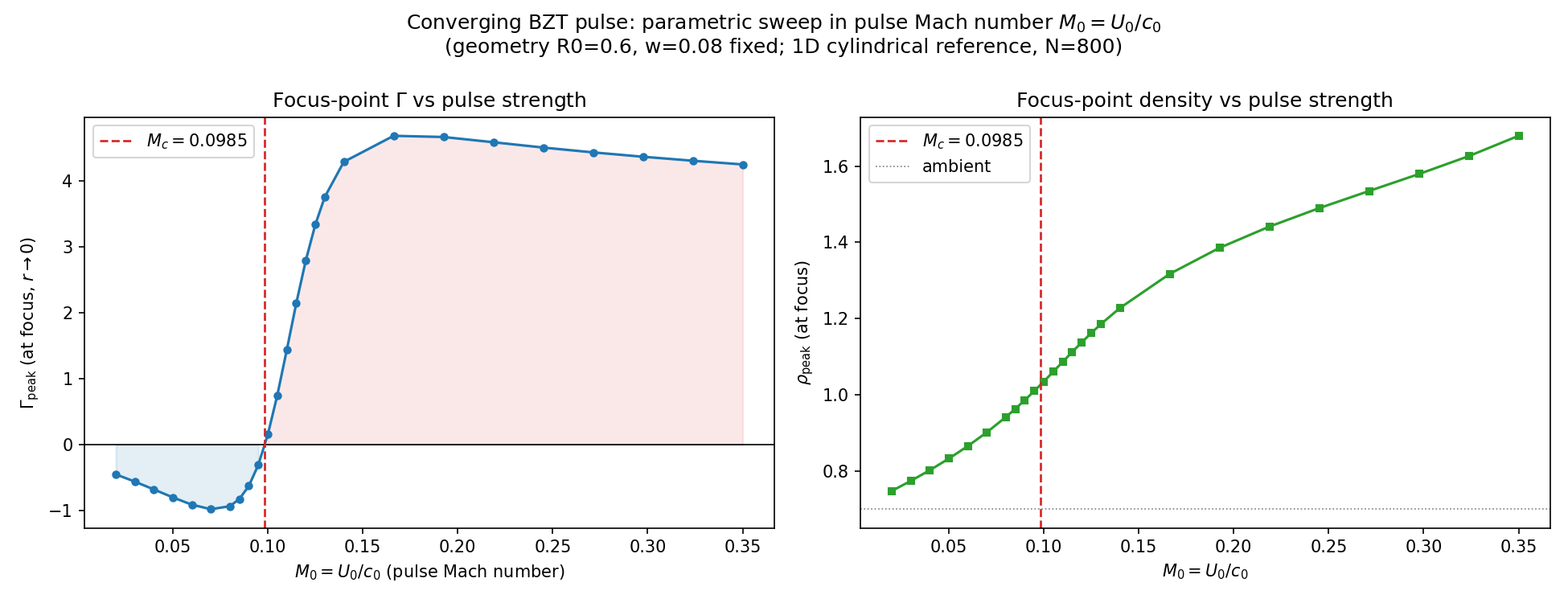}
\caption{Parametric sweep in pulse Mach number $M_0=U_0/c_{s,0}$ (1D cylindrical reference, $N=800$): focus-point $\Gamma$ (left) and focus-point density (right).}
\label{fig:bzt-mach-sweep}
\end{figure*}

The 4-quadrant problem of \S\ref{sec:val-bzt-2d} demonstrates the two signature BZT wave-type inversions (shock $\to$ compression fan, rarefaction $\to$ rarefaction shock) but, because the flow stays close to the initial isentropes, it does not probe what happens when a BZT flow is driven, by its own dynamics, from $\Gamma<0$ into the classical $\Gamma>0$ region. Geometric (cylindrical/spherical) wave focusing is a natural mechanism for exactly this: as a compression wave converges toward a symmetry axis, its amplitude grows through the purely kinematic $1/\sqrt{r}$ (cylindrical) or $1/r$ (spherical) geometric amplification, independent of the fluid's nonlinearity. If the amplification is strong enough, the local compression can push the fluid across $\Gamma=0$ even though the flow starts, and remains almost everywhere, deep in the BZT regime. This subsection designs, validates, and parametrically characterizes such a converging BZT wave, and contrasts it with the same construction in a classical ideal gas, for
which no such crossing is possible.\\
\paragraph{Governing non-dimensional number and test-case design} We first designed and cross-checked the initial condition using a cheap, independent 1D cylindrically-symmetric MUSCL(minmod)+Rusanov+RK2 solver of the geometric compressible Euler equations with a reflecting axis condition at $r=0$ ($u\to-u$) and outflow at $r=r_{\max}$. The initial condition is an axisymmetric, purely kinematic compression pulse: the van der Waals gas  ($a=1.0,b=0.3,R=1.0,\gamma=1.02$) at the uniform ambient state $\rho_0=0.70$, $P_0=0.3787$ (giving $\Gamma_0=-0.255$, deep in the BZT region, with ambient sound speed $c_{s,0}=0.4497$), on top of which an inward-pointing Gaussian velocity annulus is superposed,
\begin{equation}
\begin{split}
  u_r(r) &= -U_0\,\exp\!\left[-\left(\frac{r-R_0}{w}\right)^{2}\right], \\
  (u,v) &= u_r(r)\,\frac{(x-x_c,\,y-y_c)}{r},
\end{split}
  \label{eq:val-converging-ic}
\end{equation}
with $R_0=0.6$, $w=0.08$; the velocity vanishes at $r=0$ and in the far field, so the pulse must genuinely propagate inward and focus, rather than starting as an already-compressed core. With the geometry fixed, the single governing non-dimensional number controlling the focusing strength is the pulse Mach number
\begin{equation}
  M_0 = U_0/c_{s,0},
  \label{eq:val-converging-mach}
\end{equation}
which plays the same role here as a piston Mach number does for a planar shock: it is the only free amplitude parameter once the geometry ($R_0,w$, domain size) and the ambient thermodynamic state are fixed. Figure~\ref{fig:bzt-converging-1d-precheck} shows the 1D cylindrical solver's radial grid-converged profiles for $M_0\approx0.133$ ($U_0=0.06$), confirming the focusing event is physical. The inward annulus converges, peaks sharply near $t\approx1.4$ with $\rho$ rising from the ambient $0.70$ to $\approx1.19$ and $\Gamma$ swinging from $-0.26$ to $+3.9$ exactly at the axis, then reflects and expands outward into a rarefaction ($\rho\approx0.56$ at $r=0$, $t=1.7$) as $\Gamma$ relaxes back toward its ambient negative value.
\paragraph{2D Cartesian D2Q4 VLBM results} Having validated the mechanism in 1D, we ran the two-dimensional Cartesian counterpart with the D2Q4 VLBM, van der Waals flux on a $[0,2]\times[0,2]$ domain with the pulse of Eq.~\eqref{eq:val-converging-ic} centered at $(x_c,y_c)=(1,1)$, run at $N=800\times800$ to $t=1.7$ with outflow boundaries. Because the initial condition is exactly axisymmetric, an independent quantitative cross-check is available essentially for free: azimuthally averaging the 2D Cartesian field and comparing it against the (already grid-converged) 1D cylindrical reference at matching times, without needing to build a separate 2D reference solver.
Figure~\ref{fig:bzt-converging-2d-captures} shows the resulting density and $\Gamma$ fields near the focus. The initially uniform ambient state develops a converging annular compression ring (visible at $t=0.8$), which collapses into a sharp central peak by $t\approx1.3$--$1.4$; the $\Gamma$ field shows that the ambient blue ($\Gamma<0$) region remains essentially
undisturbed everywhere except a small, transient red ($\Gamma>0$) spot exactly at the focus at $t=1.4$ -- a genuinely two-dimensional, spatially localized realization of the same $\Gamma$-sign crossing seen in the 1D simulation. Figure~\ref{fig:bzt-converging-2d-vs-1d} confirms this quantitatively: the azimuthally-averaged 2D VLBM density profile agrees with the independent 1D cylindrical reference to within a few percent everywhere.
\paragraph{Contrast with a classical (ideal) gas} To isolate what is genuinely due to $\Gamma<0$, we reran the identical geometry, amplitude ($U_0=0.06$), and -- crucially -- the identical \emph{ambient sound speed} $c_{s,0}=0.4497$ with an ordinary ideal gas ($\gamma=1.4$, so $\Gamma\equiv(\gamma+1)/2=1.2>0$ everywhere, by
construction), by choosing the ambient pressure $P_0=c_{s,0}^2\rho_0/\gamma$.
Matching $c_{s,0}$ (rather than just $\rho_0,P_0$) means the two runs share the same linear acoustic travel time $R_0/c_{s,0}$ and the same pulse Mach number $M_0$; any difference in focusing behavior is then attributable to the EOS nonlinearity alone, not to a trivial mismatch of time scales.\\
Figure~\ref{fig:bzt-vs-classical-2d} stacks the two cases at the same absolute times. The classical run's $\Gamma$ field is, trivially, flat
at $1.2$ throughout -- no crossing is possible by construction -- but the \emph{dynamics} also differ substantially: the classical pulse peaks earlier and more weakly ($t\approx1.1$, $\rho_{\rm peak}\approx0.93$) and then overshoots directly into a strong rarefaction (visible as a
dark hole at the center by $t=1.4$--$1.5$), whereas the BZT pulse peaks later and higher ($t\approx1.4$, $\rho_{\rm peak}\approx1.10$) and stays compressed longer before its own rebound. This timing/amplitude offset, made explicit in Fig.~\ref{fig:bzt-vs-classical-history}, has a
direct thermodynamic explanation: evaluating the local sound speed along each ambient isentrope shows that in the classical gas $c_s(\rho)$
increases steadily with compression ($0.450\to0.501$ as $\rho: 0.70\to1.20$), the textbook mechanism that drives fast shock steepening, whereas in the BZT gas $c_s(\rho)$ \emph{decreases} with compression
($0.450\to0.246$ over the same range, bottoming out near $\rho\approx1.0$--$1.1$, exactly where $\Gamma$ changes sign) before
recovering once $\Gamma$ is solidly positive. A locally \emph{slower} sound speed under compression delays the characteristics
from catching up with one another, which is precisely the mechanism behind the smooth compression fans of \S\ref{sec:val-bzt-2d}, and here manifests as a later, sharper, longer-lived focus rather than the early, weak, quickly-rebounding classical peak.
\paragraph{Parametric threshold in pulse Mach number} Finally, we ran simulations to sweep the governing number $M_0=U_0/c_{s,0}$ of Eq.~\eqref{eq:val-converging-mach} at fixed geometry ($R_0=0.6,w=0.08$, $N=800$), recording the peak focus-point density $\rho_{\rm peak}$ and
fundamental derivative $\Gamma_{\rm peak}$ over $M_0\in[0.02,0.35]$ (26 values), then bisected the transition to machine precision. Figure~\ref{fig:bzt-mach-sweep} shows the result: there is a sharp critical Mach number
\begin{equation}
  M_{0,c} \approx 0.0985,
  \label{eq:val-mach-critical}
\end{equation}
below which the focus-point $\Gamma$ remains negative for all $M_0$ tested (the pulse never leaves the BZT regime, even at focus), and above
which the focus locally and transiently crosses into the classical $\Gamma>0$ region ($\Gamma_{\rm peak}$ rising from $\approx-0.3$ to
$\approx+0.75$ between $M_0=0.095$ and $M_0=0.105$). Two further features of the $\Gamma_{\rm peak}(M_0)$ curve are worth noting. First, it is non-monotonic below the crossing: $\Gamma_{\rm peak}$ initially becomes \emph{more} negative with increasing $M_0$ (down to $\approx-0.98$ near $M_0\approx0.07$) before turning around and crossing zero, i.e.\ the $\Gamma<0$ smoothing becomes transiently more effective at moderate amplitude before the geometric compression eventually overwhelms it. Second, past the crossing $\Gamma_{\rm peak}$ saturates and even declines slightly with further increasing $M_0$ (peaking $\approx4.68$ near $M_0\approx0.17$, settling to $\approx4.2$ by
$M_0=0.35$) rather than diverging -- consistent with the $\approx1.15{:}1$ pressure-ratio ceiling of this EOS's $\Gamma<0$ pocket identified when
this test case was being designed (see discussion preceding this subsection): once the focus is compressed well past the pocket, the
local thermodynamics are simply those of an ordinary dense gas, with no further BZT-specific structure to asymptote to. By contrast, the
focus-point density $\rho_{\rm peak}(M_0)$ (right panel) is smooth and monotonic straight through the crossing, with no visible feature at $M_{0,c}$ -- the transition is a genuinely thermodynamic ($\Gamma$-level) event, essentially invisible in the density field alone, which is the
central reason $\Gamma(x,y)$ rather than $\rho(x,y)$ is shown explicitly throughout this subsection and \S\ref{sec:val-bzt-2d}. The amplitude used for Figs.~\ref{fig:bzt-converging-2d-captures}--\ref{fig:bzt-vs-classical-history} ($M_0=U_0/c_{s,0}\approx0.133$) lies comfortably above $M_{0,c}$, consistent with the crossing observed there.
\subsection{Inviscid shock-vortex interaction}
\label{sec:val-shock-vortex}
\begin{figure*}[t]
\centering
\includegraphics[width=0.95\textwidth]{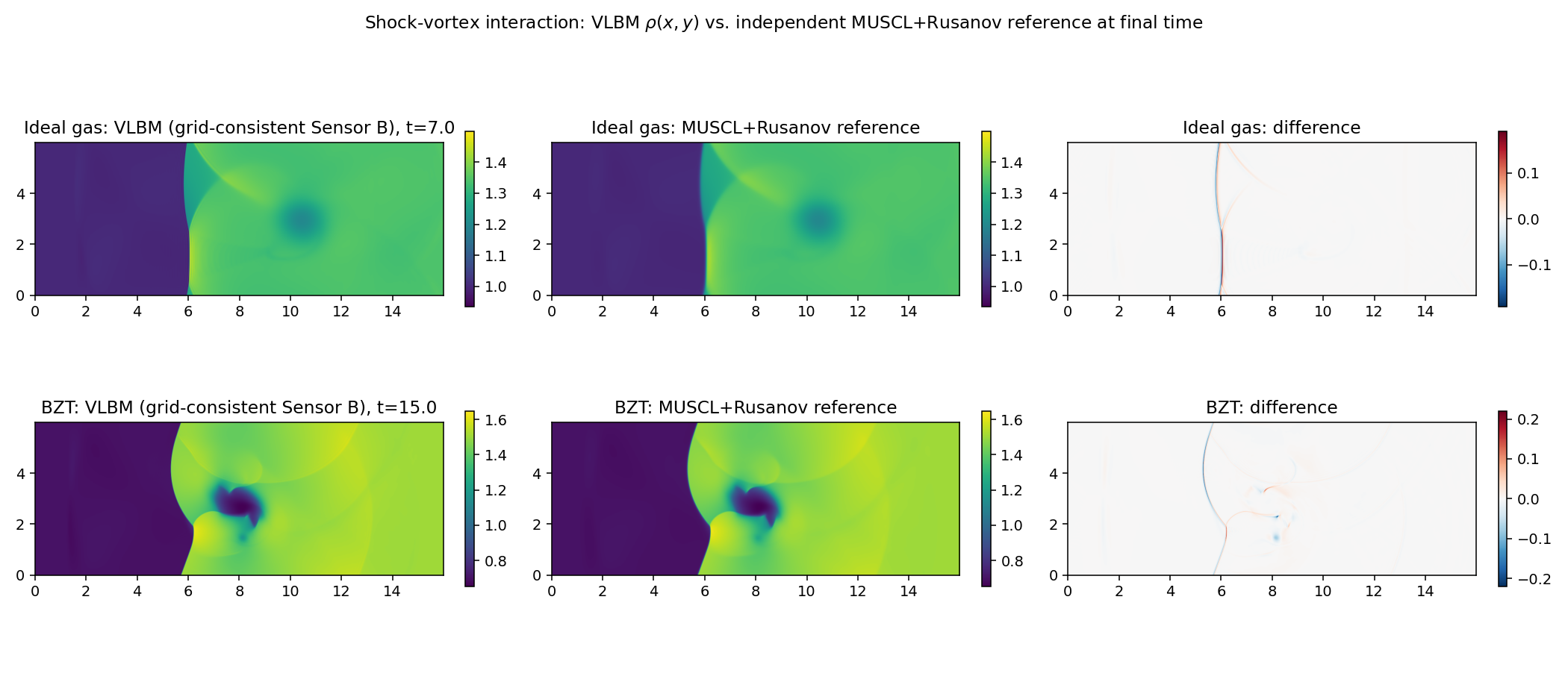}
\caption{Shock-vortex interaction, $\rho(x,y)$ at final time ($t=7$ ideal gas, top; $t=15$ BZT, bottom): D2Q4 VLBM, MUSCL+Rusanov reference, and difference.}
\label{fig:sv-vs-reference}
\end{figure*}

\begin{figure*}[t]
\centering
\includegraphics[width=0.95\textwidth]{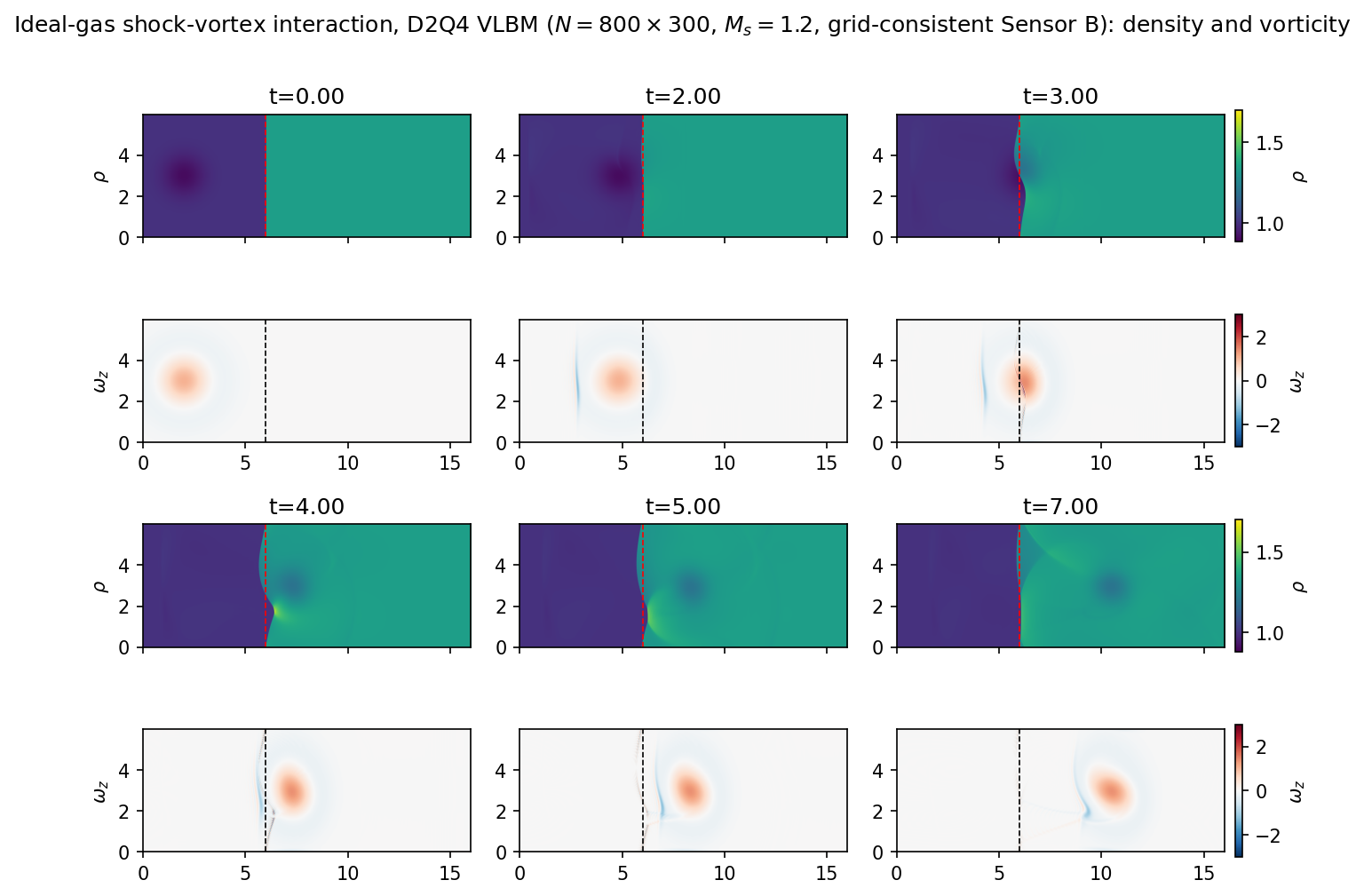}
\caption{Ideal-gas shock-vortex interaction, $M_s=1.2$, $N=800\times300$: density (top) and vorticity $\omega_z$ (bottom) at six times.}
\label{fig:sv-ideal-seq}
\end{figure*}

\begin{figure*}[t]
\centering
\includegraphics[width=0.95\textwidth]{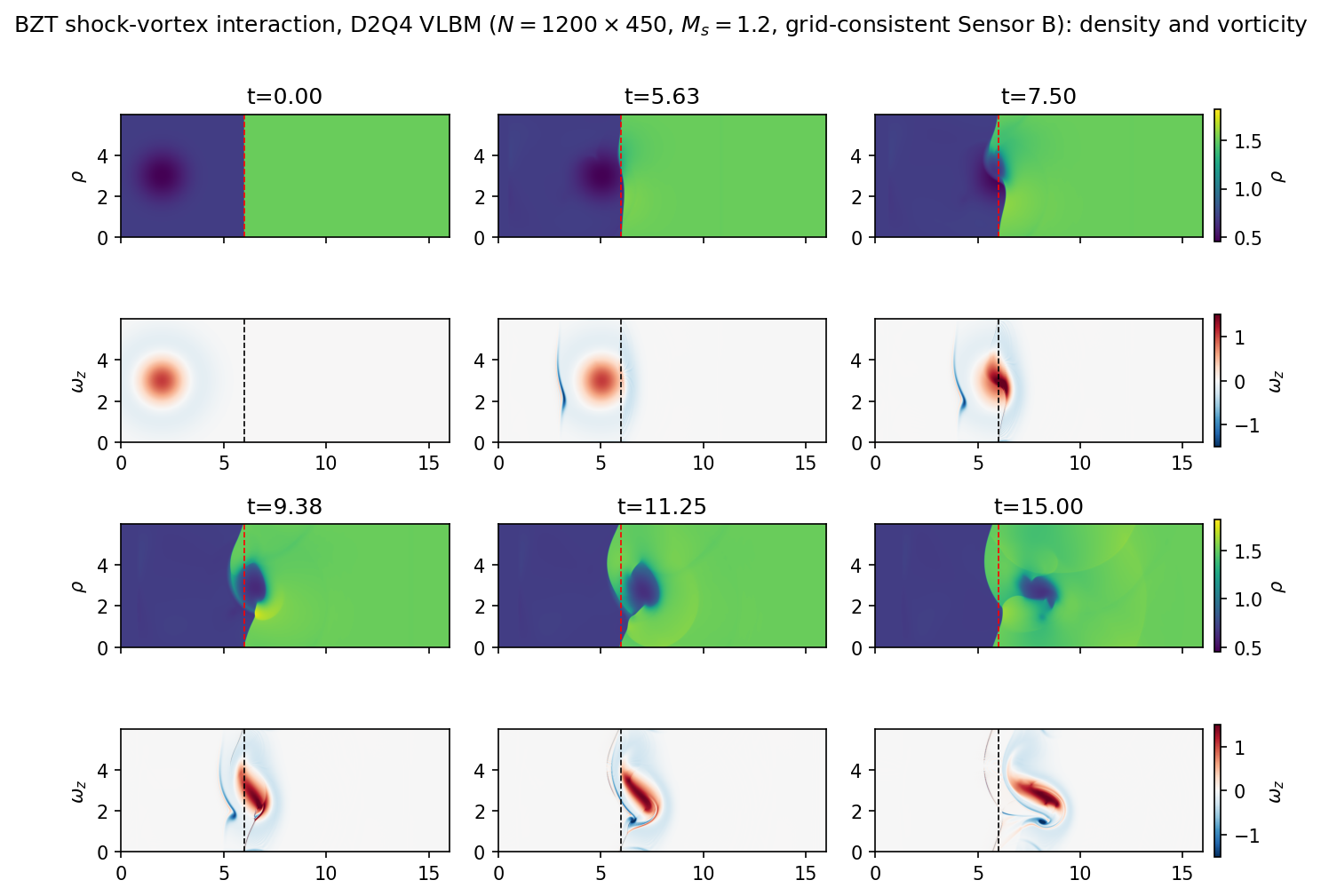}
\caption{BZT shock-vortex interaction, $M_s=1.2$, $N=1200\times450$: density (top) and vorticity $\omega_z$ (bottom) at six times.}
\label{fig:sv-bzt-seq}
\end{figure*}

The benchmarks so far have exercised shocks and smooth vortical structures separately: the two-dimensional Riemann problems of
\S\S\ref{sec:val-2d-riemann}--\ref{sec:val-bzt-2d} contain shocks and only a transient, shock-generated roll-up, while the Gresho vortex of \S\ref{sec:val-gresho} contains a smooth vortex and no shock at all. We close the validation section with a test that combines both, persistently, in the same domain: a normal shock with an isentropic vortex advecting through it. We consider identical constructions in both an ideal gas and a genuinely BZT (dense-gas, $\Gamma<0$) medium. Unlike every other benchmark in this work, we additionally cross-check the result against an independent numerical method.\\
In the shock's own rest frame, a stationary normal shock sits at $x=x_s$ on the domain $[0,L_x]\times[0,L_y]$, $L_x=16$, $L_y=6$, $x_s=6$, with upstream inflow at shock Mach number $M_s=1.2$ (left boundary: ghost cells pinned to the exact upstream Rankine-Hugoniot state; right boundary: outflow/edge-replicate; $y$: periodic). An isentropic vortex of circulation parameter $\beta_v=2$ (the conventional Yee/Shu symbol; unrelated to the relaxation parameter $\beta$ of \S\ref{sec:num-adaptive}) and core radius $R_c=1$ -- the Yee/Shu profile $u_\phi(r)=\frac{\beta_v}{2\pi R_c}\,r\,e^{\frac12(1-r^2/R_c^2)}$,
$r=|\bm x-\bm x_c|$, closed by $dP/dr=\rho u_\phi(r)^2/r$ along the ambient
isentrope -- is embedded in the upstream flow at $\bm x_c=(2,L_y/2)$ and advects into the shock. This construction is EOS-generic: for \emph{any}
EOS, a purely azimuthal $u_\phi(r)$ with $\rho(r)$ obtained by integrating $d\rho/dr=\rho u_\phi(r)^2/(r\,c_s(\rho)^2)$ inward along an isentrope from an ambient far-field value is an exact, time-independent solution of the 2D Euler equations (continuity and the azimuthal momentum equation are trivial by symmetry; only the radial momentum balance is
nontrivial, and that is exactly what is imposed). For the ideal gas this reduces to the classical closed-form Yee/Shu vortex; for the van der Waals
case it has no closed form and is instead obtained once, numerically, by RK2 shooting.\\
\emph{Ideal-gas configuration:} $\gamma=1.4$, ambient $\rho_1=1$, $P_1=1$,
giving $u_1=1.420$, $c_{s,1}=1.183$; the exact Rankine-Hugoniot state is
$\rho_2=1.342$, $u_2=1.058$, $P_2=1.513$. Run to $t_{\rm end}=7$ at
$N=800\times300$.\\
\emph{BZT configuration:} the same van der Waals gas as
\S\ref{sec:val-bzt-2d}
($a=1.0$, $b=0.3$, $R=1.0$, $\gamma=1+R/c_v=1.02$), ambient state
$\rho_1=0.70$, $T_1=0.980$ ($\approx0.993\,T_c$), giving $P_1=0.3787$,
$c_{s,1}=0.4497$, and $\Gamma_1=-0.255$ -- deep in the BZT ($\Gamma<0$) regime,
as in \S\ref{sec:val-bzt-2d}. The stationary shock is solved numerically from the exact (mass, momentum, stagnation-enthalpy) jump system at the same $M_s=1.2$: $\rho_2=1.506$, $u_2=0.251$, $P_2=0.488$, with $\Gamma_2=+4.47$
-- even this mild shock ejects the fluid into strongly classical territory,
so, unlike the isentropic compression fan of \S\ref{sec:val-bzt-2d}, the
region of genuinely $\Gamma<0$ fluid in this test is confined to the (unshocked) upstream flow and the vortex's own outer envelope, with the
vortex's centrifugal rarefaction independently pushing its core to $\Gamma>0$ from the low-density side. Because the ambient sound speed is
roughly a factor of $2.6$ lower than the ideal-gas case's, we run to a correspondingly longer $t_{\rm end}=15$, at a higher resolution,
$N=1200\times450$, to keep the vortex core comparably well resolved.\\
As throughout this work, we cross-validate against the unsplit MUSCL(minmod)+Rusanov+RK2 finite-volume scheme (EOS-generic for the BZT case, as in \S\ref{sec:val-bzt-2d}), given the same inflow/outflow/periodic boundary treatment and run at the same resolution as the VLBM ($N=800\times300$ ideal, $N=1200\times450$ BZT).
Figure~\ref{fig:sv-vs-reference} shows the final-time density field for both solutions side by side, together with the pointwise difference.
Directly differencing gives a relative $L_1$ density error of $0.093\%$
(ideal gas) and $0.114\%$ (BZT), with the residual in both cases confined to thin lines exactly on the shock front and around the vortex core -- the same pattern seen in every shock-bearing benchmark in this work.
Figures~\ref{fig:sv-ideal-seq} and~\ref{fig:sv-bzt-seq} show the density and vorticity, $\omega_z:=\mathrm{vort}(\bm x)=\partial_x u_y-\partial_y u_x$ (Eq.~\eqref{eq:num-sensor-dilvort}), fields at a sequence of times for the two configurations. In both cases the vortex advects into the shock, locally bulging it outward (the local upstream flow seen by the shock is momentarily faster where the vortex's own velocity adds to the mean flow, and slower where it subtracts), and a faint
disturbance is shed upstream into the entirely subsonic-relative-to-the-shock incoming flow -- visible in the vorticity field as thin, curved bands propagating away from the interaction point. This is the acoustic wave generated by the shock-vortex interaction, expected to be quadrupolar in the weak-vortex limit~\cite{Ellzey1995}. After crossing, the vortex re-forms downstream, visibly compressed and distorted relative to its upstream shape, in both configurations.
\section{Conclusion}
\label{sec:conclusion}
We have presented a VLBM for the compressible Euler equations with a generic, single-phase EOS. The formulation relies on the D1Q2 and D2Q4 lattices. In order to guarantee consistent convergence to the Euler limit under acoustic scaling the relaxation parameter $\beta$ is defined as a shock-capturing parameter based on gradients of the velocity, pressure and temperature fields. In addition we have introduced a consistent and conservative formulation for adaptive time-stepping. The numerical scheme has been validated using a large set of test cases in both 1D and 2D, using ideal and non-ideal EOS. Future extensions of the model, currently under development, include consistent introduction of Navier-Stokes-level dissipation and extension to adaptive grid-refinement.
\section*{Acknowledgments}
During the preparation of this work, the author(s) used Claude AI to assist with code development/data analysis and manuscript drafting. After using this tool, the author(s) reviewed and edited the content as needed and take full responsibility for the content of the published article.  This work was supported by the Swiss National Science Foundation (SNSF) Grants 200021-228065 and 200021-236715. Computational resources at the Swiss National Super Computing Center (CSCS) were provided under Grants No. s1286, sm101 and s1327.
\section*{Author Declarations}
\subsection*{Conflict of Interest}
The authors have no conflicts of interest to disclose.
\section*{Data Availability}
The data and code that support the findings of this study are available from the corresponding author upon reasonable request.
\appendix
\section{Multi-scale analysis}
\label{app:multiscale}
Starting with the discrete time-evolution equations,
\begin{equation}
  \bm f_i(\bm x+\bm c_i\delta t,\ t+\delta t) \;=\; \bm f_i(\bm x,t)
  \;+\;2\beta\Big[\bm f_i^{\rm eq}\big(\bm W(\bm x,t)\big)-\bm f_i(\bm x,t)\Big],
  \label{eq:ceg-evolution}
\end{equation}
with $\beta(\bm x,t)=\beta\big(\nabla\bm W(\bm x,t)\big)$. Taylor-expanding the shift on the left of \eqref{eq:ceg-evolution} and
inserting the ansatz,
\begin{equation}
    \bm f_i = \bm f_i^{(0)} +\delta t\,\bm f_i^{(1)} + \delta t^2\bm f_i^{(2)}+O(\delta t^3),
\end{equation}
along with,
\begin{equation}
    \frac{1}{2\beta} = \tau^{(0)}+\delta t \tau^{(1)}+\delta t^2 \tau^{(2)}+O(\delta t^3),
\end{equation}
and defining $\delta t$ as the smallness parameter for the multi-scale expansion, at order $\delta t$ we get,
\begin{equation}
    D_t^{(1)} \bm f_i^{(0)} = -\frac{1}{\tau^{(0)}} \bm f_i^{(1)}.
\end{equation}
Summing over $i$ we get the corresponding balance equations,
\begin{equation}
    \partial_t^{(1)}\bm W + \partial_\alpha \bm Q(\bm W) = 0,
\end{equation}
which are the Euler balance equations we are targeting. At the next level, i.e. $\delta t^2$ we have,
\begin{equation}
    \partial_t^{(2)}\bm f_i^{(0)} + D_t^{(1)} \bm f_i^{(1)} + \frac{1}{2} D_t^2 \bm f_i^{(0)} + \frac{\tau^{(1)}}{\tau^{(0)}}D_t \bm f_i^{(0)}= -\frac{1}{\tau^{(0)}} \bm f_i^{(2)},
\end{equation}
which using the equation from the previous order can be simplified as,
\begin{equation}
    \partial_t^{(2)}\bm f_i^{(0)} + D_t^{(1)} \left(\frac{1}{2}-\tau^{(0)}\right)D_t^{(1)}\bm f_i^{(0)} +  \frac{\tau^{(1)}}{\tau^{(0)}}D_t \bm f_i^{(0)} = -\frac{1}{\tau^{(0)}} \bm f_i^{(2)}.
\end{equation}
Summing up over $i$ we find,
\begin{equation}
  \;
  \partial_t^{(2)} \bm W
  \;-\; \sum_{\alpha=1}^{D} \partial_\alpha\!\left\{
  \left(\tau^{(0)}-\frac12\right)
  \Big[c^2\bm I - \bm A_\alpha(\bm W)^2\Big]\,\partial_\alpha\bm W
  \right\} = 0.
\end{equation}
Note that setting $\tau^{(0)}=1/2$ we get,
\begin{equation}
  \;
  \partial_t^{(2)} \bm W
  = 0.
\end{equation}
Summing up balance equations at order $\delta t$ and $\delta t^2$,
\begin{equation}
    \partial_t \bm W + \partial_\alpha \bm Q(\bm W) + O(\delta t^3) = 0.
\end{equation}
This proves that under acoustic scaling, setting $\frac{1}{2\beta}=\frac{1}{2}+\delta t\,g(\bm x,t)$ with the function $g$ guaranteed positive definite, the solver is second order convergent to the Euler limit.
\section*{references}
\bibliography{aipsamp}
\end{document}